\documentclass[onecolumn]{autart}
\usepackage{graphicx,color}
\usepackage{amsmath,amssymb,mathrsfs,amsfonts,theorem}
\usepackage{algorithm,algorithmic}
\renewcommand{\algorithmicensure}{\textbf{Input:}}
\usepackage{subfigure}

\newcommand{\Ker}{\operatorname{Ker}}
\newcommand{\Image}{\operatorname{Im}}
\newcommand{\Basis}{\operatorname{Basis}}

\newcommand{\E}{\operatorname{\mathbb{E}}}

\newcommand{\V}{\mathcal{V}}

\newtheorem{theorem}{Theorem}[section]
\newtheorem{lemma}{Lemma}[section]
\newtheorem{corollary}{Corollary}[section]
\newtheorem{remark}{Remark}

\newtheorem{problem}{Problem}

\renewcommand{\natural}{{\mathbb{N}}}
\renewcommand{\int}{{\mathbb{Z}}}
\newcommand{\real}{{\mathbb{R}}}

\newcommand{\transpose}{\mathsf{T}} %or \top or \intercal
\DeclareMathOperator*{\argmin}{arg\,min}

\newcommand{\until}[1]{\{1,\dots,#1\}}
\newcommand{\Supp}{S}
\newcommand{\Dec}{D}

\renewcommand{\theenumi}{(\roman{enumi})}

\begin{document}

\begin{frontmatter}

\title{Distributed Estimation via Iterative Projections 
with Application to Power Network Monitoring}

  % \title{A Distributed Projection Method for Power
  %   Network Estimation and Monitoring}

  % \title{A Distributed Method for State Estimation and False Data
  %   Detection in Power Networks}

  % \title{Distributed Estimation and False Data Detection\\ with
  %   Application to Power Networks}
  
  \thanks{This material is based in part upon work supported by ICB ARO
    grant W911NF-09-D-0001 and NSF grant CNS-0834446, and in part upon the
    EC Contract IST 224428 "CHAT".}

\author[First]{Fabio Pasqualetti}\ead{fabiopas@engineering.ucsb.edu},
\author[Second]{Ruggero Carli}\ead{carlirug@dei.unipd.it},
%\author[Second]{Antonio Bicchi}\ead{bicchi@ing.unipi.it},
\author[First]{Francesco Bullo}\ead{bullo@engineering.ucsb.edu}

\address[First]{\mbox{Center} for \mbox{Control}, \mbox{Dynamical}
  \mbox{Systems} and \mbox{Computation}, University of California,
  Santa Barbara, USA \\\tt
  \{fabiopas,bullo\}@engineering.ucsb.edu}

\address[Second]{Departement of Information Engineering, University of
  Padova, Padova, Italy \\\tt carlirug@dei.unipd.it}

\begin{abstract}
  This work presents a distributed method for control centers to
  monitor the operating condition of a power network, i.e., to
  estimate the network state, and to ultimately determine the
  occurrence of threatening situations. State estimation has been
  recognized to be a fundamental task for network control centers to
  ensure correct and safe functionalities of power grids. We consider
  (static) state estimation problems, in which the state vector
  consists of the voltage magnitude and angle at all network buses. We
  consider the state to be linearly related to network measurements,
  which include power flows, current injections, and voltages phasors
  at some buses. We admit the presence of several cooperating control
  centers, and we design two distributed methods for them to compute
  the minimum variance estimate of the state given the network
  measurements. The two distributed methods rely on different modes of
  cooperation among control centers: in the first method an
  \emph{incremental} mode of cooperation is used, whereas, in the
  second method, a \emph{diffusive} interaction is implemented. Our
  procedures, which require each control center to know only the
  measurements and structure of a subpart of the whole network, are
  computationally efficient and scalable with respect to the network
  dimension, provided that the number of control centers also
  increases with the network cardinality. Additionally, a
  \emph{finite-memory} approximation of our diffusive algorithm is
  proposed, and its accuracy is characterized. Finally, our estimation
  methods are exploited to develop a distributed algorithm to detect
  corrupted data among the network measurements.
\end{abstract}

\end{frontmatter}

\section{Introduction}\label{sec:introduction}
% \begin{enumerate}
%   \item application to false data detection and power networks
%   \item include monitoring
%   \item change a little abstract to match with title, and leave the
%     fact that we first do estimation, and then have a suitable
%     application
%   \item in 1.2, reference also Kaczmarz methods (see paragraph in
%     section 1)
%   \item find adjective for estimation in title (row action? ART?
%     projection?)
% \end{enumerate}
Large-scale complex systems, such as, for instance, the electrical power
grid and the telecommunication system, are receiving increasing attention
from researchers in different fields. The wide spatial distribution and the
high dimensionality of these systems preclude the use of centralized
solutions to tackle classical estimation, control, and fault detection
problems, and they require, instead, the development of new decentralized
techniques. One possibility to overcome these issues is to geographically
deploy some monitors in the network, each one responsible for a different
subpart of the whole system. Local estimation and control schemes can
successively be used, together with an information exchange mechanism to
recover the performance of a centralized scheme. 
% The focus of this work is on distributed estimation schemes for
% power systems.  As it will be clear in the next sections, however,
% the methods we develop are general and have applicability beyond the
% considered scenario (see \cite{FB-RC-AB-FB:10n} for an application
% to linear dynamical networks).

\subsection{Control centers, state estimation and cyber security in power
  networks}
\begin{figure}[tb]
  \centering \subfigure[]{
    \includegraphics[width=.45\columnwidth]{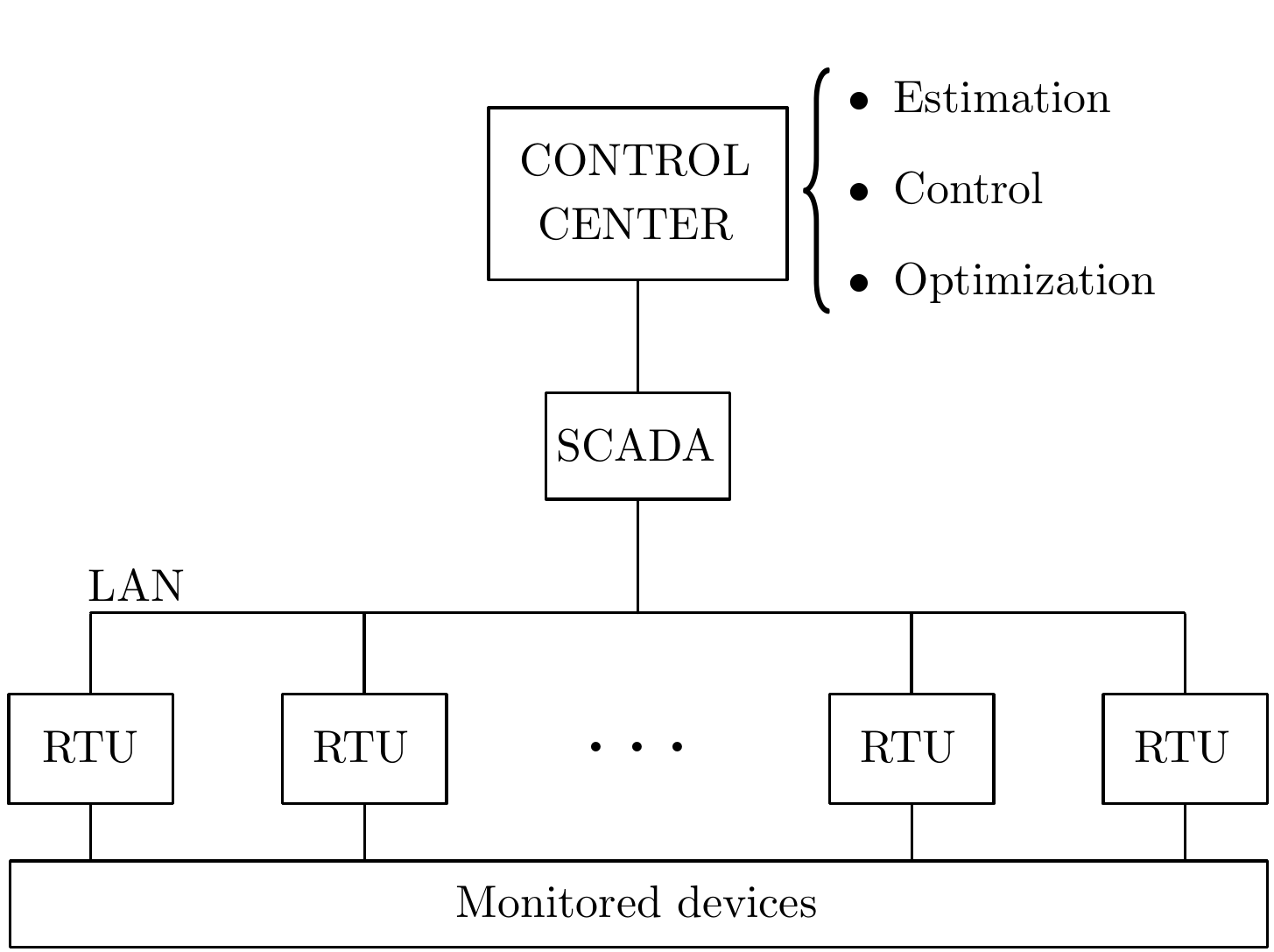}
    \label{fig:control_center}
  } \subfigure[]{
    \includegraphics[width=.49\columnwidth]{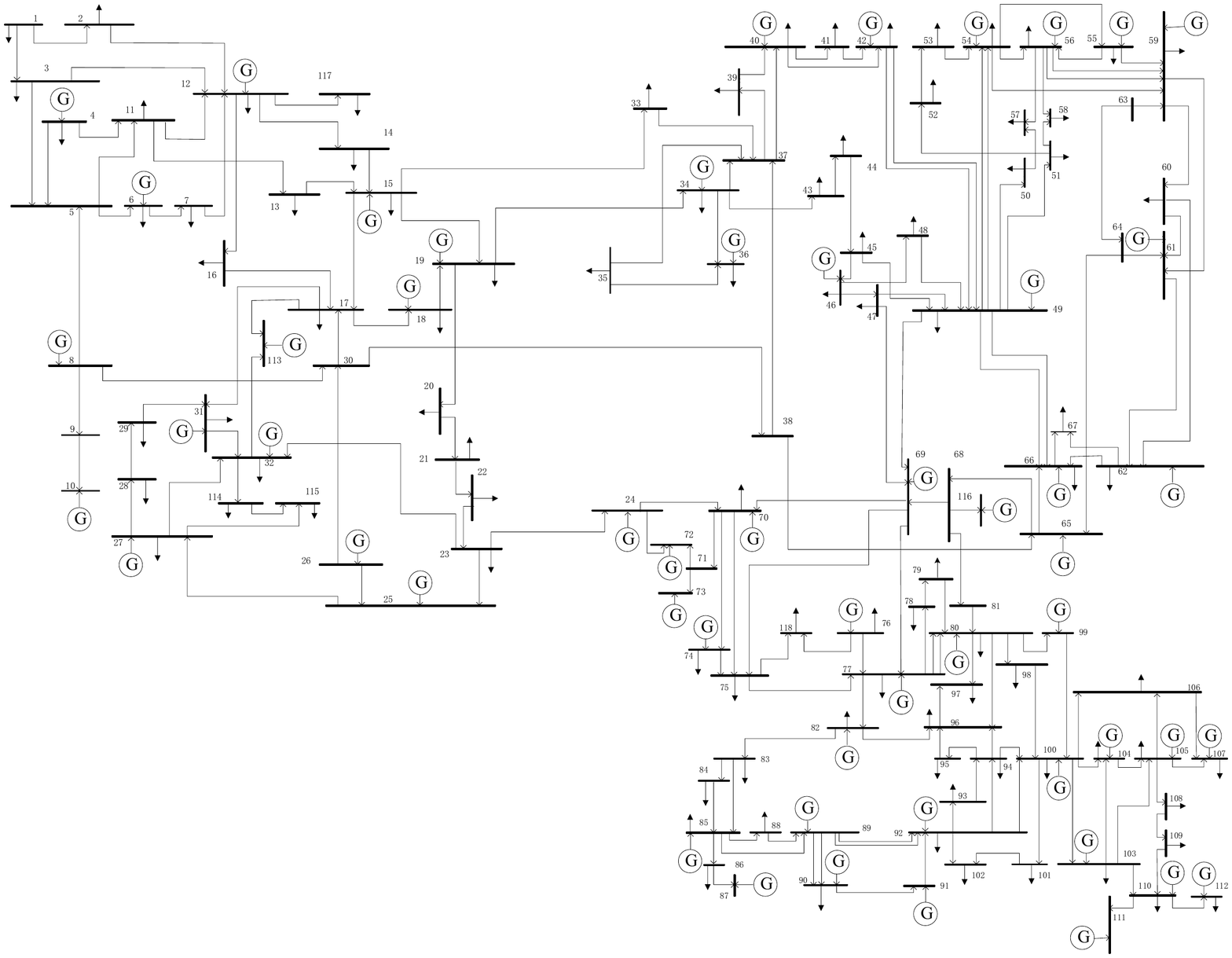}
    \label{fig:power_network118}
  }
  \caption[Optional caption for list of figures]{In
    Fig.~\ref{fig:control_center}, remote terminal units (RTU) transmit
    their measurements to a SCADA terminal via a LAN network. The data is
    then sent to a Control Center to implement network estimation, control,
    and optimization procedures. Fig.~\ref{fig:power_network118} shows the
    diagram of the IEEE 118 bus system (courtesy of the IIT Power
    Group). The network has 118 buses, 186 branches, 99 loads, and 54
    generators.}
\end{figure}

Power systems are operated by system operators from the area control
center. The main goal of the system operator is to maintain the network in
a secure operating condition, in which all the loads are supplied power by
the generators without violating the operational limits on the transmission
lines. In order to accomplish this goal, at a given point in time, the
network model and the phasor voltages at every system bus need to be
determined, and preventive actions have to be taken if the system is found
in an insecure state. For the determination of the operating state, remote
terminal units and measuring devices are deployed in the network to gather
measurements. These devices are then connected via a local area network to
a SCADA (Supervisory Control and Data Acquisition) terminal, which supports
the communication of the collected measurements to a control center. At the
control center, the measurement data is used for control and optimization
functions, such as contingency analysis, automatic generation control, load
forecasting, optimal power flow computation, and reactive power dispatch
\cite{AA-AGS:04}. A diagram representing the interconnections between
remote terminal units and the control center is reported in
Fig.~\ref{fig:control_center}.
% \begin{figure}
%     \centering
%     \includegraphics[width=.4\columnwidth]{./control_center}
%     \caption{Remote terminal units (RTU) transmit their measurements
%       to a SCADA terminal via a LAN network. The data is then sent to
%       a Control Center to implement network estimation, control, and
%       optimization procedures.}
%     \label{fig:control_center}
% \end{figure}
Various sources of uncertainties, e.g., measurement and communication
noise, lead to inaccuracies in the received data, which may affect the
performance of the control and optimization algorithms, and,
ultimately, the stability of the power plant. This concern was first
recognized and addressed in
\cite{FCS-JW:70,FCS-JW:70-bis,FCS-JW:70-bisbis} by introducing the
idea of (static) state estimation in power systems.

Power network state estimators are broadly used to obtain an optimal
estimate from redundant noisy measurements, and to estimate the state
of a network branch which, for economical or computational reasons, is
not directly monitored. For the power system state estimation problem,
several centralized and parallel solutions have been developed in the
last decades, e.g., see \cite{AM:99,DMF-FFW-LM:95,MS-YW:03}. Being an
online function, computational issues, storage requirements, and
numerical robustness of the solution algorithm need to be taken into
account. Within this regard, distributed algorithms based on network
partitioning techniques are to be preferred over centralized ones. %
Moreover, even in decentralized setting, the work in \cite{N:04} on
the blackout of August 2003 suggests that an estimation of the entire
network is essential to prevent networks damages. In other words, the
whole state vector should be estimated by and available to every
unit. %
The references \cite{LZ-AA:05,WJ-VV-GTH:07} explore the idea of using
a global control center to coordinate estimates obtained locally by
several local control centers. In this work, we improve upon these
prior results by proposing a fully decentralized and distributed
estimation algorithm, which, by only assuming local knowledge of the
network structure by the local control centers, allows them to obtain
in finite time an optimal estimate of the network state. Being the
computation distributed among the control centers, our procedure
appears scalable against the power network dimension, and,
furthermore, numerically reliable and accurate.

A second focus of this paper is false data detection and cyber attacks
in power systems.  Because of the increasing reliance of modern power
systems on communication networks, the possibility of cyber attacks is
a real threat \cite{JM:07}. One possibility for the attacker is to
corrupt the data coming from the measuring units and directed to the
control center, in order to introduce arbitrary errors in the
estimated state, and, consequently, to compromise the performance of
control and optimization algorithms \cite{YL-MKR-PN:09}. This
important type of attack is often referred in the power systems
literature to as \emph{false data injection attack}. Recently, the
authors of \cite{LX-YM-BS:10} show that a false data injection attack,
in addition to destabilizing the grid, may also lead to fluctuations
in the electricity market, causing significant economical losses. The
presence of false data is classically checked by analyzing the
statistical properties of the \emph{estimation residual} $z - H \hat
x$, where $z$ is the measurements vector, $\hat x$ is a state
estimate, and $H$ is the state to measurements matrix. For an attack
to be successful, the residual needs to remain within a certain
confidence level. Accordingly, one approach to circumvent false data
injection attacks is to increase the number of measurements so as to
obtain a more accurate confidence bound. Clearly, by increasing the
number of measurements, the data to be transmitted to the control
center increases, and the dimension of the estimation problem
grows. By means of our estimation method, we address this
dimensionality problem by distributing the false data detection
problem among several control centers.

\subsection{Related work on distributed estimation and projection
  methods}
Starting from the eighties, the problem of distributed estimation has
attracted intense attention from the scientific community, generating
along the years a very rich literature. More recently, because of the
advent of highly integrated and low-cost wireless devices as key
components of large autonomous networks, the interest for this
classical topic has been renewed. For a wireless sensor network, novel
applications requiring efficient distributed estimation procedures
include, for instance, environment monitoring, surveillance,
localization, and target tracking. Considerable effort has been
devoted to the development of distributed and adaptive filtering
schemes, which generalize the notion of adaptive estimation to a setup
involving networked sensing and processing devices
\cite{FSC-CGL-AHS:08}. In this context, relevant methods include
incremental Least Mean-Square \cite{CGL-AHS:07}, incremental Recursive
Least-Square \cite{AHS-CGL:07}, Diffusive Least Mean-Square
\cite{AHS-CGL:07}, and Diffusive Recursive
Least-Square\cite{FSC-CGL-AHS:08}. Diffusion Kalman filtering and
smoothing algorithms are proposed, for instance, in
\cite{RC-AC-LS-SZ:08a,FSC-AHS:10}, and consensus based techniques in
\cite{IS-AR-GG:07,IDS-GM-GBG:09}. We remark that the strategies
proposed in the aforementioned references could be adapted for the
solution of the power network static estimation problem. Their
assumptions, however, appear to be not well suited in our context for
the following reasons. First, the convergence of the above estimation
algorithms is only asymptotic, and it depends upon the communication
topology. As a matter of fact, for many communication topologies, such
as Cayley graphs and random geometric graphs, the convergence rate is
very slow and scales badly with the network dimension. Such slow
convergence rate is clearly undesirable because a delayed state
estimation could lead the power plant to instability. Second,
approaches based on Kalman filtering require the knowledge of the
global state and observation model by all the components of the
network, and they violate therefore our assumptions. Third and
finally, the application of these methods to the detection of cyber
attacks, which is also our goal, is not straightforward, especially
when detection guarantees are required. An exception is constituted by
\cite{SSS-MSS-DMS:09}, where a estimation technique based on local
Kalman filters and a consensus strategy is developed. This latter
method, however, besides exhibiting asymptotic convergence, does not
offer guarantees on the final estimation error.

Our estimation technique belongs to the family of Kaczmarz
(row-projection) methods for the solution of a linear system of
equations. See \cite{SK:37,RG-RB-GTH:70,KT:71,YC:81} for a detailed
discussion. Differently from the existing row-action methods, our
algorithms exhibit finite time convergence towards the exact solution,
and they can be used to compute any weighted least squares solution to
a system of linear equations.

\subsection{Our contributions}
The contributions of this work are threefold. First, we adopt the
static state network estimation model, in which the state vector is
linearly related to the network measurements. We develop two methods
for a group of interconnected control centers to compute an optimal
estimate of the system state via distributed computation. Our first
estimation algorithm assumes an \emph{incremental} mode of cooperation
among the control centers, while our second estimation algorithm is
based upon a \emph{diffusive} strategy. Both methods are shown to
converge in a finite number of iterations, and to require only local
information for their implementation. Differently than
\cite{CR-SP-SU:05}, our estimation procedures assume neither the
measurement error covariance nor the measurements matrix to be
diagonal. Furthermore, our algorithms are advantageous from a
communication perspective, since they reduce the distance between
remote terminal units and the associated control center, and from a
computational perspective, since they distribute the measurements to
be processed among the control centers. Second, as a minor
contribution, we describe a finite-time algorithm to detect via
distributed computation if the measurements have been corrupted by a
malignant agent. Our detection method is based upon our state
estimation technique, and it inherits its convergence
properties. Notice that, since we assume the measurements to be
corrupted by noise, the possibility exists for an attacker to
compromise the network measurements while remaining undetected (by
injecting for instance a vector with the same noise statistics). With
respect to this limitation, we characterize the class of corrupted
vectors that are guaranteed to be detected by our procedure, and we
show optimality with respect to a centralized detection
algorithm. Third, we study the scalability of our methods in networks
of increasing dimension, and we derive a finite-memory approximation
of our diffusive estimation strategy. For this approximation procedure
we show that, under a reasonable set of assumptions and independent of
the network dimension, each control center is able to recover a good
approximation of the state of a certain subnetwork through little
computation. Moreover, we provide bounds on the approximation error
for each subnetwork. Finally, we illustrate the effectiveness of our
procedures on the IEEE 118 bus system.

The rest of the paper is organized as follows. In Section
\ref{sec:setup} we introduce the problem under consideration, and we
describe the mathematical setup. Section \ref{sec:solver} contains our
main results on the state estimation and on the detection problem, as
well as our algorithms. Section \ref{sec:approximate} describes our
approximated state estimation algorithm. In Section
\ref{sec:numerical} we study the IEEE 118 bus system, and we present
some simulation results. Finally, Section \ref{sec:conclusion}
contains our conclusion.

\section{Problem setup and preliminary notions}\label{sec:setup}
% \begin{figure}
%     \centering
%     \includegraphics[width=.6\columnwidth]{./power_network118}
%     \caption{Courtesy of the IIT Power Group. Diagram of the IEEE 118
%       bus system. The network has 118 buses, 186 branches, 99 loads,
%       and 54 generators.}
%     \label{fig:power_network118}
% \end{figure}
For a power network, an example of which is reported in
Fig. \ref{fig:power_network118}, the state at a certain instant of
time consists of the voltage angles and magnitudes at all the system
buses. The (static) state estimation problem introduced in the seminal
work by Schweppe \cite{FCS-JW:70} refers to the procedure of
estimating the state of a power network given a set of measurements of
the network variables, such as, for instance, voltages, currents, and
power flows along the transmission lines. To be more precise, let $x
\in \real^n$ and $z \in \real^p$ be, respectively, the state and
measurements vector. Then, the vectors $x$ and $z$ are related by the
relation
\begin{align}\label{nonlinear_estimation}
  z = h(x) + \eta,
\end{align}
where $h(\cdot)$ is a nonlinear measurement function, and where
$\eta$, which is traditionally assumed to be a zero mean random vector
satisfying $\E [\eta \eta^\transpose] = \Sigma_{\eta} =
\Sigma_{\eta}^\transpose >0$, is the noise measurement. An optimal
estimate of the network state coincides with the \emph{most likely}
vector $\hat x$ that solves equation \eqref{nonlinear_estimation}. It
should be observed that, instead of by solving the above estimation
problem, the network state could be obtained by measuring directly the
voltage phasors by means of phasor measurement
devices.\footnote{Phasor measurement units are devices that
  synchronize by using GPS signals, and that allow for a direct
  measurement of voltage and current phasors.} Such an approach,
however, would be economically expensive, since it requires to deploy
a phasor measurement device at each network bus, and it would be very
vulnerable to communication failures \cite{AA-AGS:04}. In this work,
we adopt the approximated estimation model presented in
\cite{FCS-JW:70-bis}, which follows from the linearization around the
origin of equation \eqref{nonlinear_estimation}. Specifically,
\begin{align}\label{linear_estimation}
  z = H x + v,
\end{align}
where $H \in \real^{p \times n}$ and where $v$, the noise measurement,
is such that $\E[v] = 0$ and$E[ v v^\transpose ] = \Sigma =
\Sigma^\transpose >0$.
% {\color{blue} In particular,
% letting $H \in \real^{p \times n}$, $\E[v] = 0$, and $E[ v
% v^\transpose ] = \Sigma = \Sigma^\transpose >0$, the state estimation
% problem becomes
% \begin{align}\label{linear_estimation}
%   z = H x + v.
% \end{align}
%}
Observe that, because of the interconnection structure of a power
network, the measurement matrix $H$ is usually sparse. Let $\Ker(H)$
denote the null space defined by the matrix $H$. For the equation
\eqref{linear_estimation}, without affecting generality, assume
$\Ker(H) = \{0\}$, and recall from \cite{DGL:69} that the vector
\begin{align}\label{eq:minimum_estimate}
  x_{\textup{wls}} = (H^\transpose \Sigma^{-1} H)^{-1} H^\transpose \Sigma^{-1} z
\end{align}
minimizes the weighted variance of the estimation
error, i.e., $x_{\textup{wls}} = \argmin_{\hat x} (z - H \hat
x)^\transpose \Sigma^{-1} (z - H \hat x)$.

The centralized computation of the minimum variance estimate to
\eqref{linear_estimation} assumes the complete knowledge of the
matrices $H$ and $\Sigma$, and it requires the inversion of the matrix
$H^\transpose \Sigma^{-1} H$. For a large power network, such
computation imposes a limitation on the dimension of the matrix $H$,
and hence on the number of measurements that can be efficiently
processed to obtain a real-time state estimate. Since the performance
of network control and optimization algorithms depend upon the
precision of the state estimate, a limitation on the network
measurements constitutes a bottleneck toward the development of a more
efficient power grid. A possible solution to address this complexity
problem is to distribute the computation of $x_{\textup{wls}}$ among
geographically deployed control centers (monitors), in a way that each
monitor is responsible for a subpart of the whole network. To be more
precise, let the matrices $H$ and $\Sigma$, and the vector $z$ be
partitioned as\footnote{In most application the error covariance
  matrix is assumed to be diagonal, so that each submatrix $\Sigma_i$
  is very sparse. However, we do not impose any particular structure
  on the error covariance matrix.}
\begin{align}\label{matrix_decompose}
  H=\left[
    \begin{array}{c}
      H_1\\
      H_2\\
      \vdots\\
      H_m
    \end{array}
  \right],\quad
  \Sigma=\left[
    \begin{array}{c}
      \Sigma_1\\
      \Sigma_2\\
      \vdots\\
      \Sigma_m
    \end{array}
  \right],\quad
  z=\left[
    \begin{array}{c}
      z_1\\
      z_2\\
      \vdots\\
      z_m
    \end{array}
  \right],
\end{align}
where, for $i \in \{ 1,\dots,m \}$, $m_i \in \natural$, $H_i \in
\real^{m_i \times n}$, $\Sigma_i \in \real^{m_i \times p}$, $z_i \in
\real^{m_i}$, and $\sum_{i=1}^m m_i = p$. Let $G=(V,\mathcal{E})$ be a
connected graph in which each vertex $i \in V = \{ 1,\dots, m\}$
denotes a monitor, and $\mathcal{E} \in V \times V$ denotes the set of
monitors interconnections. For $i\in \until{m}$, assume that monitor
$i$ knows the matrices $H_i$, $\Sigma_i$, and the vector
$z_i$. Moreover, assume that two neighboring monitors are allowed to
cooperate by exchanging information. Notice that, if the full matrices
$H$ and $\Sigma$ are nowhere available, and if they cannot be used for
the computation of $x_{\textup{wls}}$, then, with no cooperation among
the monitors, the vector $x_{\textup{wls}}$ cannot be computed by any
of the monitor. Hence we consider the following problem.

\begin{problem}[Distributed state estimation]
  Design an algorithm for the monitors to compute the minimum variance
  estimate of the network state via distributed computation.
\end{problem}

We now introduce the second problem addressed in this work. Given the
distributed nature of a power system and the increasing reliance on
local area networks to transmit data to a control center, there exists
the possibility for an attacker to compromise the network
functionality by corrupting the measurements vector. When a malignant
agent corrupts some of the measurements, the state to measurements
relation becomes
\begin{align*}
   z = H x + v + w,
\end{align*}
where the vector $w \in \real^p$ is chosen by the attacker, and,
therefore, it is unknown and unmeasurable by any of the monitoring
stations. We refer to the vector $w$ to as \emph{false} data. From the
above equation, it should be observed that there exist vectors $w$
that cannot be detected through the measurements $z$. For instance, if
the false data vector is intentionally chosen such that $w \in
\Image(H)$, then the attack cannot be detected through the
measurements $z$. Indeed, denoting with $^\dag$ the pseudoinverse
operation, the vector $x + H^\dag w$ is a valid network state. In this
work, we assume that the vector $w$ is detectable from the
measurements $z$, and we consider the following problem.

\begin{problem}[Distributed detection]
  Design an algorithm for the monitors to detect the presence of
  false data in the measurements via distributed computation.
\end{problem}

As it will be clear in the sequel, the complexity of our methods
depends upon the dimension of the state, as well as the number of
monitors. In particular, few monitors should be used in the absence of
severe computation and communication contraints, while many monitors
are preferred otherwise. We believe that a suitable choice of the
number of monitors depends upon the specific scenario, and it is not
further discussed in this work.

\begin{remark}[Generality of our methods]
  In this paper we focus on the state estimation and the false data
  detection problem for power systems, because this field of research
  is currently receiving sensible attention from different
  communities. The methods described in the following sections,
  however, are general, and they have applicability beyond the power
  network scenario. For instance, our procedures can be used for state
  estimation and false data detection in dynamical system, as
  described in \cite{FB-RC-AB-FB:10n} for the case of sensors
  networks.
\end{remark}

\section{Optimal state estimation and false data detection via
  distributed computation}\label{sec:solver}
The objective of this section is the design of distributed methods to
compute an optimal state estimate from measurements. With respect to a
centralized method, in which a powerful central processor is in charge
of processing all the data, our procedures require the computing units
to have access to only a subset of the measurements, and are shown to
reduce significantly the computational burden. In addition to being
convenient for the implementation, our methods are also optimal, in
the sense that they maintain the same estimation accuracy of a
centralized method.

For a distributed method to be implemented, the interaction structure
among the computing units needs to be defined. Here we consider two
modes of cooperations among the computing units, and, namely, the
\emph{incremental} and the \emph{diffusive} interaction. In an
incremental mode of cooperation, information flows in a sequential
manner from one node to the adjacent one. This setting, which usually
requires the least amount of communications \cite{MGR-RDN:05}, induces
a cyclic interaction graph among the processors. In a diffusive
strategy, instead, each node exchanges information with all (or a
subset of) its neighbors as defined by an interaction graph. In this
case, the amount of communication and computation is higher than in
the incremental case, but each node possesses a good estimate before
the termination of the algorithm, since it improves its estimate at
each communication round. This section is divided into three parts. In
Section \ref{subsec:solver}, we first develop a distributed
incremental method to compute the minimum norm solution to a set of
linear equations, and then exploit such method to solve a minimum
variance estimation problem. In Section \ref{sec:asynchronous} we
derive a diffusive strategy which is amenable to asynchronous
implementation. Finally, in Section \ref{sec:detection} we propose a
distributed algorithm for the detection of false data among the
measurements. Our detection procedure requires the computation of the
minimum variance state estimate, for which either the incremental or
the diffusive strategy can be used.

\subsection{Incremental solution to a set of linear equations}
We start by introducing a distributed incremental procedure to compute
the minimum norm solution to a set of linear equations. This procedure
constitutes the key ingredient of the incremental method we later
propose to solve the minimum variance estimation problem.
\begin{algorithm}[t!]
  \caption{\textit{Incremental minimum norm solution (i-th monitor)}}
  \label{algo:pseudoinverse}
  \begin{algorithmic}
   \ENSURE $H_i$, $z_i$;
    
   \REQUIRE $[z_1^\transpose \, \dots \, z_m^\transpose]^\transpose
   \in \Image(  [H_1^\transpose \, \dots \, H_m^\transpose]^\transpose )$;
    
   \IF{$i = 1$}
   
   \STATE $\hat x_0 :=0$, $K_0 := I_n$;

   \ELSE

   \STATE receive $\hat x_{i-1}$ and $K_{i-1}$ from monitor $i-1$;

   \ENDIF
   
   \smallskip
   
   \STATE $\hat x_i := \hat x_{i-1} + K_{i-1} (H_i K_{i-1})^\dag (z_i
   - H_i \hat x_{i-1})$;
    
   \smallskip

   \STATE $K_{i} := \Basis( K_{i-1}\Ker(H_i K_{i-1}) )$;
   
   \smallskip

   \IF{$i < m$}

   \STATE transmit $\hat x_i$ and $K_{i}$ to monitor $i+1$;

   \ELSE

   \RETURN $\hat x_m$;

   \ENDIF

 \end{algorithmic}
\end{algorithm}
Let $H \in \real^{p\times m}$, and let $z \in \Image(H)$, where
$\Image(H)$ denotes the range space spanned by the matrix
$H$. Consider the system of linear equations $z = H x$, and recall
that the unique minimum norm solution to $z = H x$ coincides with the
vector $\hat x$ such that $z=H \hat x$ and $\|\hat x\|_2$ is
minimum. It can be shown that $\|\hat x\|_2$ being minimum corresponds
to $\hat x$ being orthogonal to the null space $\Ker(H)$ of $H$
\cite{DGL:69}. Let $H$ and $z$ be partitioned in $m$ blocks as in
\eqref{matrix_decompose}, and let $G=(V,\mathcal{E})$ be a directed
graph such that $V = \until{m}$ corresponds to the set of monitors,
and, denoting with $(i,j)$ the directed edge from $j$ to $i$,
$\mathcal{E} = \{(i+1,i) : i = 1,\dots,m-1\} \cup \{(1,m)\}$. Our
incremental procedure to compute the minimum norm solution to $z=H
\hat x$ is in Algorithm \ref{algo:pseudoinverse}, where, given a
subspace $\V$, we write $\Basis(\V)$ to denote any full rank matrix
whose columns span the subspace $\V$.
 % Assume that the following operations are performed by the $i$-th
 % monitor:
 % \begin{enumerate}
 % \item receive the vector $\hat x_{i-1}$ and the matrix $K_{i-1}$ from
 %   monitor $i-1$
 % \item compute the vector $\hat x_i = \hat x_{i-1} + K_{i-1} (H_i
 %   K_{i-1})^\dag (z_i - H_i \hat x_{i-1})$,
 %   % \begin{align*}
 %   %   \hat x_i = \hat x_{i-1} + K_{i-1} (H_i K_{i-1})^\dag (z_i - H_i
 %   %   \hat x_{i-1}),
 %   % \end{align*}
 % \item compute the matrix $K_i$ such that $\Image(K_i) = K_{i-1} \Ker(H_i K_{i-1})$,
 %   % \begin{align*}
 %   %   \Image(K_i) = K_{i-1} \Ker(H_i K_{i-1}),
 %   % \end{align*}
 % \item transmit $\hat x_i$ and $K_i$ to monitor $i+1$.
 % \end{enumerate}
We now proceed with the analysis of the convergence properties of the
\emph{Incremental minimum norm solution} algorithm.

\begin{theorem}\emph{(Convergence of Algorithm
    \ref{algo:pseudoinverse})}\label{thm:Algo:Pseudo}
  Let $z=Hx$, where $H$ and $z$ are partitioned in $m$ row-blocks as
  in \eqref{matrix_decompose}. In Algorithm \ref{algo:pseudoinverse},
  the $m$-th monitor returns the vector $\hat x$ such that $z=H\hat x$
  and $\hat x \perp \Ker(H)$.
\end{theorem}
\begin{pf}
  See Section \ref{pf_algopseudo}.
\end{pf}

It should be observed that the dimension of $K_i$ decreases, in
general, when the index $i$ increases. In particular, $K_m = \{ 0 \}$
and $K_1 = \Ker(H_1)$. To reduce the computational burden of the
algorithm, monitor $i$ could transmit the smallest among
$\Basis(K_{i-1} \Ker(H_i K_{i-1}))$ and $\Basis(K_{i-1} \Ker(H_i
K_{i-1})^\perp)$, together with a packet containing the type of the
transmitted basis.

\begin{remark}\emph{\bf (Computational complexity of Algorithm
  \ref{algo:pseudoinverse})}
In Algorithm \ref{algo:pseudoinverse}, the main operation to be
performed by the $i$-th agent is a singular value decomposition
(SVD).\footnote{The matrix $H$ is usually very sparse, since it
  reflects the network interconnection structure. Efficient SVD
  algorithms for very large sparse matrices are being developed
  (cf. \emph{SVDPACK}).} Indeed, since the range space and the null
space of a matrix can be obtained through its SVD, both the matrices
$(H_i K_{i-1})^\dag$ and $\Basis ( K_{i-1} \Ker (H_i K_{i-1}) )$ can
be recovered from the SVD of $H_i K_{i-1}$. Let $H \in \real^{m \times
  n}$, $m > n$, and assume the presence of $\lceil m / k \rceil$
monitors, $1 \le k \le m $. Recall that, for a matrix $M \in \real^{k
  \times p}$, the singular value decomposition can be performed with
complexity $O ( \textup{min} \{k p^2 , k^2 p \} )$
\cite{GHG-CFvL:89}. Hence, the computational complexity of computing a
minimum norm solution to the system $z = Hx$ is $O(m n^2)$. In Table
\ref{table:complexity} we report the computational complexity of
Algorithm \ref{algo:pseudoinverse} as a function of the size $k$.

\begin{table}[ht]
  \caption{Computational complexity of Algorithm \ref{algo:pseudoinverse}.}
  \centering
  \bigskip
  \begin{tabular}{cccc}
    \hline\hline
    Block size & $I$-th complexity & Total complexity & Communications\\
    \hline
    $k \le n$ & $O(k^2 n)$ & $O(m k n)$ & $\lceil m / k \rceil-1$\\
    $k > n$ & $O(k n^2)$ & $O(m n^2)$ & $\lceil m / k \rceil-1$\\
    \hline
  \end{tabular}
  \label{table:complexity}
  \bigskip
\end{table}

The following observations are in order. First, if $k \le n$, then the
computational complexity sustained by the $i$-th monitor is much
smaller than the complexity of a centralized implementation, i.e., $O
( k^2 n ) \ll O( m n^2 )$. Second, the complexity of the entire
algorithm is optimal, since, in the worst case, it maintains the
computational complexity of a centralized solution, i.e., $O ( m k n )
\le O(m n ^2)$. Third and finally, a compromise exists between the
blocks size $k$ and the number of communications needed to terminate
Algorithm \ref{algo:pseudoinverse}. In particular, if $k = m$, then no
communication is needed, while, if $k = 1$, then $m-1$ communication
rounds are necessary to terminate the estimation
algorithm.\footnote{Additional $m-1$ communication rounds are needed
  to transmit the estimation to every other monitor.}
\end{remark}

\subsection{Incremental state estimation via distributed
  computation}\label{subsec:solver}
We now focus on the computation of the weighted least squares solution
to a set of linear equations. Let $v$ be an unknown and unmeasurable
random vector, with $\E(v)=0$ and $\E(vv^\transpose) = \Sigma =
\Sigma^\transpose > 0$. Consider the system of equations
\begin{align}\label{system_noise}
  z=Hx + v,
\end{align}
and assume $\Ker(H)=0$. Notice that, because of the noise vector $v$,
we generally have $z \not\in \Image(H)$, so that Algorithm
\ref{algo:pseudoinverse} cannot be directly employed to compute the
vector $x_{\textup{wls}}$ defined in \eqref{eq:minimum_estimate}. It
is possible, however, to recast the above weighted least squares
estimation problem to be solvable with Algorithm
\ref{algo:pseudoinverse}. Note that, because the matrix $\Sigma$ is
symmetric and positive definite, there exists\footnote{Choose for
  instance $B = W\Lambda^{1/2}$, where $W$ is a basis of eigenvectors
  of $\Sigma$ and $\Lambda$ is the corresponding diagonal matrix of
  the eigenvalues.} a full row rank matrix $B$ such that
$\Sigma=BB^\transpose$. Then, equation \eqref{system_noise} can be
rewritten as
\begin{align}\label{system_under}
  z=\left[
    \begin{array}{cc}
      H & \varepsilon B
    \end{array}
    \right]
    \left[
    \begin{array}{c}
      x \\
      \bar v
    \end{array}
    \right],
\end{align}
where $\varepsilon \in \real_ {>0}$, $\E[\bar v]=0$ and $\E[\bar v \bar
v^\transpose]=\varepsilon^{-2} I$. Observe that, because $B$ has full
row rank, the system \eqref{system_under} is underdetermined, i.e., $z
\in \Image([H \enspace \varepsilon B])$ and $\Ker([H \enspace
\varepsilon B]) \neq 0$. Let
\begin{align}\label{eq:x_epsilon}
  \left[
    \begin{array}{c}
      \hat x(\varepsilon)\\
      \hat{\bar v}
    \end{array}
  \right]=
  \left[
    \begin{array}{cc}
      H & \varepsilon B
    \end{array}
    \right]^\dag z.
\end{align}
The following theorem characterizes the relation between the minimum
variance estimation $x_\textup{wls}$ and $\hat x(\varepsilon)$.
\begin{theorem}\emph{(Convergence with $\varepsilon$)}\label{existence_limit}
  Consider the system of linear equations $z = Hx + v$. Let $\E(v)=0$
  and $\E(vv^\transpose)=\Sigma = BB^\transpose > 0$ for a full row
  rank matrix $B$. Let
\begin{align*}
  \begin{split}
    C&=\varepsilon (I - HH^\dag) B ,\qquad
    E=I-C^\dag C,\qquad
    D=\varepsilon E [ I + \varepsilon^2 E B^\transpose
    (HH^\transpose)^\dag B E ]^{-1} B^\transpose (HH^\transpose)^\dag
    (I -\varepsilon B C^\dag).
    \end{split}
\end{align*}

Then
\begin{align*}
  \left[
    \begin{array}{cc}
      H & \varepsilon B
    \end{array}
  \right]^\dag = \left[
    \begin{array}{c}
      H^\dag - \varepsilon H^\dag B (C^\dag + D)  \\
      C^\dag + D
    \end{array}
  \right];
\end{align*}
and
\begin{align*}
  \lim_{\varepsilon \rightarrow 0^+} H^\dag - \varepsilon H^\dag B
  (C^\dag + D) = (H^\transpose \Sigma^{-1}
  H)^{-1}H^\transpose\Sigma^{-1}.
\end{align*}
\end{theorem}
\begin{pf}
  See Section \ref{proof_existence_limit}.
\end{pf}

Throughout the paper, let $\hat x(\varepsilon)$ be the vector defined
in \eqref{eq:x_epsilon}, and notice that Theorem \ref{existence_limit}
implies that
\begin{align*}
  x_{\textup{wls}}=\lim_{\varepsilon \rightarrow 0^+} \hat x(\varepsilon).
\end{align*}

\begin{remark}[Incremental state estimation]\label{remark:incremental}
  For the system of equations $z=Hx+v$, let $BB^\transpose$ be the
  covariance matrix of the noise vector $v$, and let
    \begin{align}\label{incremental_partition}
      H=\left[
    \begin{array}{c}
      H_1\\
      H_2\\
      \vdots\\
      H_m
    \end{array}
  \right],\quad
   B=\left[
    \begin{array}{c}
      B_1\\
      B_2\\
      \vdots\\
      B_m
    \end{array}
  \right],\quad
  z=\left[
    \begin{array}{c}
      z_1\\
      z_2\\
      \vdots\\
      z_m
    \end{array}
  \right],
  \end{align}
  where $m_i \in \natural$, $H_i \in \real^{m_i \times n}$, $B_i \in
  \real^{m_i \times p}$, and $z_i \in \real^{m_i}$. For $\varepsilon >
  0$, the estimate $\hat x(\varepsilon)$ of the weighted least squares
  solution to $z=Hx+v$ can be computed by means of Algorithm
  \ref{algo:pseudoinverse} with input $[H_i \; \varepsilon B_i]$ and
  $z_i$.
  % the estimate $\hat
  % x(\varepsilon)$ can be computed by means of Algorithm
  % \ref{algo:pseudoinverse} with input $[H_i \; \varepsilon B_i]$ and
  % $z_i$, $i \in \until{m}$, where, being $BB^\transpose$ is the noise
  % covariance matrix, we have
  % \begin{align}\label{incremental_partition}
  %     H=\left[
  %   \begin{array}{c}
  %     H_1\\
  %     H_2\\
  %     \vdots\\
  %     H_m
  %   \end{array}
  % \right],\quad
  %  B=\left[
  %   \begin{array}{c}
  %     B_1\\
  %     B_2\\
  %     \vdots\\
  %     B_m
  %   \end{array}
  % \right],\quad
  % z=\left[
  %   \begin{array}{c}
  %     z_1\\
  %     z_2\\
  %     \vdots\\
  %     z_m
  %   \end{array}
  % \right],
  % \end{align}
  % and, $m_i \in \natural$, $H_i \in \real^{m_i \times n}$, $B_i \in
  % \real^{m_i \times n}$, $z_i \in \real^{m_i}$, and $\varepsilon > 0$.
\end{remark}

Observe now that the estimate $\hat x(\varepsilon)$ coincides with
$\hat x_{\textup{wls}}$ only in the limit for $\varepsilon \rightarrow
0^+$. When the parameter $\varepsilon$ is fixed, the estimate $\hat
x(\varepsilon)$ differs from the minimum variance estimate $\hat
x_{\textup{wls}}$. We next characterize the approximation error
$x_{\textup{wls}} - \hat x(\varepsilon)$.

\begin{corollary}\emph{(Approximation error)}\label{approx_error}
  Consider the system $z=Hx+v$, and let
  $\E[vv^\transpose]=BB^\transpose$ for a full row rank matrix
  $B$. Then
  \begin{align*}
    x_{\textup{wls}} - \hat x(\varepsilon) = \varepsilon H^\dag B D z,
  \end{align*}
  where $D$ is as in Theorem \ref{existence_limit}.
\end{corollary}
\begin{pf}
  With the same notation as in the proof of Theorem
  \ref{existence_limit}, for every value of $\varepsilon > 0$, the
  difference $x_{\textup{wls}} - \hat x(\varepsilon)$ equals
%  \begin{align*}
    $\left((H^\transpose \Sigma^{-1}
    H)^{-1} H^\transpose \Sigma^{-1} - H^\dag + \varepsilon H^\dag B
    (C^\dag+D)  \right)z$.
 % \end{align*}
  Since $(H^\transpose \Sigma^{-1} H)^{-1} H^\transpose \Sigma^{-1} -
  H^\dag + \varepsilon H^\dag B C^\dag = 0$ for every $\varepsilon >
  0$, it follows $x_{\textup{wls}} - \hat x(\varepsilon) = \varepsilon
  H^\dag B D z$.
\end{pf}

Therefore, for the solution of system \eqref{system_noise} by means of
Algorithm \ref{algo:pseudoinverse}, the parameter $\varepsilon$ is
chosen according to Corollary \ref{approx_error} to meet a desired
estimation accuracy. It should be observed that, even if the entire
matrix $H$ needs to be known for the computation of the exact
parameter $\varepsilon$, the advantages of our estimation technique
are preserved. Indeed, if the matrix $H$ is unknown and an upper bound
for $\|H^\dag BDz\|$ is known, then a value for $\varepsilon$ can
still be computed that guarantees the desired estimation accuracy. On
the other hand, even if $H$ is entirely known, it may be inefficient
to use $H$ to perform a centralized state estimation over
time. Instead, the parameter $\varepsilon$ needs to be computed only
once. To conclude this section, we characterize the estimation
residual \mbox{$z -H \hat x$}. This quantity plays an important role
for the synthesis of a distributed false data detection algorithm.

\begin{corollary}\emph{(Estimation residual)}\label{residual}
  Consider the system $z=Hx+v$, and let $\E[vv^\transpose] = \Sigma =
  \Sigma^\transpose > 0$. Then\footnote{Given a vector $v$ and a
    matrix $H$, we denote by $\|v\|$ any vector norm, and by $\|H\|$
    the corresponding induced matrix norm.}
  \begin{align*}
    \lim_{\varepsilon \rightarrow 0^+} \| z - H \hat x (\varepsilon)
    \| \le \| (I - H W)\|\| v \|,
  \end{align*}
  where $W = (H^\transpose\Sigma^{-1}H)^{-1}H^\transpose\Sigma^{-1}$.
\end{corollary}
\begin{pf}
  By virtue of Theorem \ref{existence_limit} we have
  %\begin{align*}
    $\lim_{\varepsilon \rightarrow 0^+} \hat x (\varepsilon) =
    x_{\textup{wls}} =
    (H^\transpose\Sigma^{-1}H)^{-1}H^\transpose\Sigma^{-1} z = W z$.
  %\end{align*}
  Observe that $H WH=H$, and recall that $z = Hx + v$. For any matrix
  norm, we have
  \begin{align*}
    \|z - H x_{\textup{wls}}\| = \|z - H W z\|
    = \| (I - H W) (H x + v) \|
    =\| H x - H x +(I - H W) v \|
    \le \| (I - H W)\| \| v \|,
  \end{align*}
  and the theorem follows.
\end{pf}

\subsection{Diffusive state estimation via distributed
  computation}\label{sec:asynchronous}
The implementation of the incremental state estimation algorithm described
in Section \ref{subsec:solver} requires a certain degree of coordination
among the control centers. For instance, an ordering of the monitors is
necessary, such that the $i$-th monitor transmits its estimate to the
$(i+1)$-th monitor. This requirement imposes a constraint on the monitors
interconnection structure, which may be undesirable, and, potentially, less
robust to link failures. In this section, we overcome this limitation by
presenting a diffusive implementation of Algorithm
\ref{algo:pseudoinverse}, which only requires the monitors interconnection
structure to be connected.\footnote{An undirected graph is said to be
  connected if there exists a path between any two vertices
  \cite{CDG-GFR:01}.} To be more precise, let $V = \{1,\dots,m\}$ be the
set of monitors, and let $G = (V,E)$ be the undirected graph describing the
monitors interconnection structure, where $E \subseteq V \times V$, and
$(i,j) \in E$ if and only if the monitors $i$ and $j$ are connected. The
neighbor set of node $i$ is defined as $N_i = \{j \in V : (i,j) \in
E\}$. We assume that $G$ is connected, and we let the distance between two
vertices be the minimum number of edges in a path connecting them. Finally,
the diameter of a graph $G$, in short $\textup{diam}(G)$, equals the
greatest distance between any pair of vertices. Our diffusive procedure is
described in Algorithm \ref{algo:finite_solver}, where the matrices $H_i$
and $\varepsilon B_i$ are as defined in equation
\eqref{incremental_partition}.  During the $h$-th iteration of the
algorithm, monitor $i$, with $i\in \until{N}$, performs the following three
actions in order:
\begin{enumerate}
\item transmits its current estimates $\hat{x}_i$ and $K_i$ to all
  its neighbors;
\item receives the estimates $\hat{x}_j$ from neighbors $N_i$; and
\item updates $\hat{x}_i$ and $K_i$ as in the \emph{for} loop of
  Algorithm \ref{algo:finite_solver}.
\end{enumerate}

We next show the convergence of Algorithm
\ref{algo:finite_solver} to the minimum variance estimate.

\begin{algorithm}[t!]
  \caption{\textit{Diffusive state estimation ($i$-th monitor)}}
  \label{algo:finite_solver}
  \begin{algorithmic}
    
    \ENSURE $H_i$, $\varepsilon B_i$, $z_i$;

    %\STATE $n:=$ number of columns of $H_i$;

    \STATE $\hat x_i := [H_i \; \varepsilon B_i]^\dag z_i$;

    \STATE $K_i := \Basis(\Ker([H_i \; \varepsilon B_i]))$;

    \WHILE{$K_i \neq 0$}
    
    \FOR{$j \in N_i$}
    
    \STATE receive $\hat x_j$ and $K_j$;

    \STATE $\hat x_i := \hat x_i + [K_i \enspace 0] [-K_i \enspace
    K_j]^\dag (\hat x_i - \hat x_j)$;

    \STATE $K_i := \Basis(\Image(K_i) \cap \Image(K_j))$;
    
    \ENDFOR
    
    \STATE transmit $\hat x_i$ and $K_i$;

    \ENDWHILE

  \end{algorithmic}
\end{algorithm}

% \begin{algorithm}[t]
%   \SetKwInOut{Input}{Input}
%   \SetKwInOut{Require}{Require}
%   \SetKwInOut{Set}{Define}
%   \SetKwInOut{Title}{Algorithm}
%   \Input{$H_i$, $\varepsilon B_i$, $z_i$\;}
%   \Set{$n:=$ number of columns of $H_i$\;}

%   $\hat x_i := [H_i \; \varepsilon B_i]^\dag z_i$\;
  
%   $K_i := \Basis(\Ker([H_i \; \varepsilon B_i]))$\;

%   \While{$K_i \neq 0$}{
      
%     \For{$j \in N_i$}{

%       receive $\hat x_j$ and $K_j$\;

%       $\hat x_i := \hat x_i + [K_i \enspace 0] [-K_i \enspace K_j]^\dag (\hat x_i -
%       \hat x_j)$\;

%       $K_i := \Basis(\Image(K_i) \cap \Image(K_j))$\;
%     } 
%     transmit $\hat x_i$ and $K_i$\;
%   }
%   \Return{$\hat x_i(1:n)$\;}
%   \caption{\textit{Asynchronous state estimation ($i$-th monitor)}}
%   \label{algo:finite_solver}
% \end{algorithm}

\begin{theorem}[Convergence of Algorithm \ref{algo:finite_solver}]
  Consider the system of linear equations $z = Hx + v$, where $\E [v]
  = 0$ and $\E[v v^\transpose] = B B^\transpose$. Let $H$, $B$ and $z$
  be partitioned as in \eqref{incremental_partition}, and let
  $\varepsilon > 0$. Let the monitors communication graph be
  connected, let $d$ be its diameter, and let the monitors execute the
  Diffusive state estimation algorithm. Then, each monitor computes
  the estimate $\hat x(\varepsilon)$ of $x$ in $d$ steps.
\end{theorem}
\begin{pf}
  Let $\hat x_i$ be the estimate of the monitor $i$, and let $K_i$ be
  such that $x - \hat x_i \in \Image(K_i)$, where $x$ denotes the
  network state, and $\hat x_i \perp \Image(K_i)$. Notice that $z_i =
  [H_i \; \varepsilon B_i] \hat x_i$, where $z_i$ it the $i$-th
  measurements vector. Let $i$ and $j$ be two neighboring
  monitors. Notice that there exist vectors $v_i$ and $v_j$ such that
  $\hat x_i + K_i v_i = \hat x_j + K_j v_j$. In particular, those
  vectors can be chosen as
  \begin{align*}
    \left[
      \begin{array}{c}
        v_i\\
        v_j
      \end{array}
    \right]=[-K_i \; K_j]^\dag (\hat x_i - \hat x_j).
  \end{align*}
  It follows that the vector 
  \begin{align*}
    \hat x_i^+ = \hat x_i + [K_i \enspace 0] [-K_i \enspace K_j]^\dag
    (\hat x_i - \hat x_j)
  \end{align*}
  is such that $z_i = [H_i \; \varepsilon B_i] \hat x_i^+$ and $z_j =
  [H_j \; \varepsilon B_j] \hat x_i^+$.  Moreover we have $\hat x_i^+
  \perp (\Image(K_i) \cap \Image(K_j))$. Indeed, notice that
    \begin{align*}
      \begin{bmatrix}
        v_i \\ v_j
      \end{bmatrix}
      \perp \Ker([-K_i \; K_j])
      \supseteq
      \left\{
        \begin{bmatrix}
          w_i \\ w_j
        \end{bmatrix}
        \,:\, K_i w_i = K_j w_j \right\}.
    \end{align*}
    % let
    % \begin{align*}
    %   \begin{bmatrix}
    %     w_i \\ w_j
    %   \end{bmatrix}
    %   =
    %   [-K_i \; K_j]^\dag (\hat x_i - \hat x_j) \perp \Ker([-K_i \; K_j]),
    % \end{align*}
    % and notice that
    % \begin{align*}
    %   \Ker\left(
    %   \begin{bmatrix}
    %     -K_i & K_j
    %   \end{bmatrix}
    %   \right)
    %   \supseteq
    %   \left\{
    %     \begin{bmatrix}
    %       v_i \\ v_j
    %     \end{bmatrix}
    %     \,:\, K_i v_i = K_j v_j \right\}.
    % \end{align*}
    We now show that $K_i v_i \perp \Image(K_j)$. By contradiction, if
    $K_i v_i \not\perp \Image(K_j)$, then $v_i = \tilde v_i +\bar
    v_i$, with $K_i \tilde v_i \perp \Image(K_j)$ and $K_i \bar v_i
    \in \Image(K_j)$. Let $\bar v_j = K_j^\dag K_i \bar v_i$, and
    $\tilde v_j = v_j - \bar v_j$. Then, $[\bar v_i^\transpose \; \bar
    v_j^\transpose]^\transpose \in \Ker([-K_i \; K_j])$, and hence
    $[v_i^\transpose \; v_j^\transpose]^\transpose \not\perp
    \Ker([-K_i \; K_j])$, which contradicts the hypothesis. We
    conclude that $[K_i \enspace 0] [-K_i \enspace K_j]^\dag (\hat x_i
    - \hat x_j) \perp \Image(K_j)$, and, since $\hat x_i \perp
    \Image(K_i)$, it follows $\hat x_i^+ \perp (\Image(K_i) \cap
    \Image(K_j))$. The theorem follows from the fact that after a
    number of steps equal to the diameter of the monitors
    communication graph, each vector $\hat x_i$ verifies all the
    measurements, and $\hat x_i \perp \Image(K_1) \cap \dots \cap
    \Image(K_m)$.
% Moreover, we have $\hat x_i^+ \perp \Ker([-K_i \; K_j])$, and
% consequently $\hat x_i^+ \perp (\Image(K_i) \cap
% \Image(K_j))$. Indeed, $\hat x_i \perp \Image(K_i)$, and, since
% $(\Image(K_i) \cap \Image(K_j)) \subseteq \Ker([-K_i \; K_j])$, it
% follows $[K_i \enspace 0] [-K_i \enspace K_j]^\dag (\hat x_i - \hat
% x_j) \perp \Image(K_j)$. The theorem follows from the fact that after
% a number of steps equal to the diameter of the monitors communication
% graph, each vector $\hat x_i$ verifies all the measurements, and $\hat
% x_i \perp \Image(K_1) \cap \dots \cap \Image(K_m)$.
\end{pf}

As a consequence of Theorem \ref{existence_limit}, in the limit for
$\varepsilon$ to zero, Algorithm \ref{algo:finite_solver} returns the
minimum variance estimate of the state vector, being therefore the
diffusive counterpart of Algorithm \ref{algo:pseudoinverse}. A detailed
comparison between incremental and diffusive methods is beyond the purpose
of this work, and we refer the interested reader to
\cite{CGL-AHS:07,CGL-AHS:08} and the references therein for a thorough
discussion. Here we only underline some key differences. While Algorithm
\ref{algo:pseudoinverse} requires less operations, being therefore
computationally more efficient, Algorithm \ref{algo:finite_solver} does not
constraint the monitors communication graph. Additionally, Algorithm
\ref{algo:finite_solver} can be implemented adopting general asynchronous
communication protocols. For instance, consider the \emph{Asynchronous
  (diffusive) state estimation} algorithm, where, at any given instant of
time in $\natural$, at most one monitor, say $j$, sends its current
estimates to its neighbors, and where, for $i\in N_j$, monitor $i$ performs
the following operations:
  \begin{enumerate}
  \item $\hat x_i := \hat x_i + [K_i \enspace 0] [-K_i \enspace
    K_j]^\dag (\hat x_i - \hat x_j)$,
    \item $K_i := \Basis(\Image(K_i) \cap \Image(K_j))$.
  \end{enumerate}
\begin{corollary}[Asynchronous estimation]
  Consider the system of linear equations $z = Hx + v$, where $\E [v] = 0$
  and $\E[v v^\transpose] = B B^\transpose$. Let $H$, $B$ and $z$ be
  partitioned as in \eqref{incremental_partition}, and let $\varepsilon >
  0$. Let the monitors communication graph be connected, let $d$ be its
  diameter, and let the monitors execute the Asynchronous (diffusive) state
  estimation algorithm. Assume that there exists a duration $T\in \natural$
  such that, within each time interval of duration $T$, each monitor
  transmits its current estimates to its neighbors. Then, each monitor
  computes the estimate $\hat x(\varepsilon)$ of $x$ within time $dT$.
\end{corollary}
\begin{pf}
  The proof follows from the following two facts. First, the
  intersection of subspaces is a commutative operation. Second, since
  each monitor performs a data transmission within any time interval
  of length $T$, it follows that, at time $dT$, the information
  related to one monitor has propagated through the network to every
  other monitor.
\end{pf}

\subsection{Detection of false data via distributed
  computation}\label{sec:detection}
In the previous sections we have shown how to compute an optimal state
estimate via distributed computation. A rather straightforward
application of the proposed state estimation technique is the
detection of false data among the measurements. When the measurements
are corrupted, the state to measurements relation becomes
\begin{align*}
  z = H x + v + w,
\end{align*}
where $w$ is the false data vector. As a consequence of Corollary
\ref{residual}, the vector $w$ is detectable if it affects
significantly the estimation residual, i.e., if
%\begin{align*}
$\lim_{\varepsilon \rightarrow 0} \| z - H \hat x(\varepsilon) \| > \Gamma$,
%\end{align*}
where the threshold $\Gamma$ depends upon the magnitude of the noise
$v$. Notice that, because false data can be injected at any time by a
malignant agent, the detection algorithm needs to be executed over
time by the control centers. Let $z(t)$ be the measurements vector at
a given time instant $t$, and let $\E[z(t_1)z^\transpose (t_2)] = 0$
for all $t_1 \neq t_2$. Based on this considerations, our distributed
detection procedure is in Algorithm \ref{algo:detection}, where the
matrices $H_i$ and $\varepsilon B_i$ are as defined in equation
\eqref{incremental_partition}, and $\Gamma$ is a predefined threshold.

\begin{algorithm}[t!]
  \caption{\textit{False data detection ($i$-th monitor)}}
  \label{algo:detection}
  \begin{algorithmic}

    \ENSURE $H_i$, $\varepsilon B_i$, $\Gamma$;

    \WHILE{True}

    \STATE collect measurements $z_i(t)$;

    \STATE estimate network state $\hat x(t)$ via Algorithm
    \ref{algo:pseudoinverse} or \ref{algo:finite_solver}\;

    \IF{$\| z_i(t) - H_i \hat x(t) \|_{\infty} > \Gamma$}

    \RETURN False data detected;

    \ENDIF

    \ENDWHILE

  \end{algorithmic}
\end{algorithm}

In Algorithm \ref{algo:detection}, the value of the threshold $\Gamma$
determines the false alarm and the misdetection rate. Clearly, if
$\Gamma \ge \| (I - H W)\|\| v(t) \|$ and $\varepsilon$ is
sufficiently small, then no false alarm is triggered, at the expenses
of the misdetection rate. By decreasing the value of $\Gamma$ the
sensitivity to failures increases together with the false alarm
rate. Notice that, if the magnitude of the noise signals is bounded by
$\gamma$, then a reasonable choice of the threshold is $\Gamma =
\gamma \| (I - H W)\|_{\infty}$, where the use of the infinity norm in
Algorithm \ref{algo:detection} is also convenient for the
implementation. Indeed, since the condition $\|z(t) - H \hat
x(t)\|_{\infty}>\Gamma$ is equivalent to $\|z_i(t) - H_i \hat
x(t)\|_{\infty}>\Gamma$ for some monitor $i$, the presence of false
data can be independently checked by each monitor without further
computation. Notice that an eventual alarm message needs to be
propagated to all other control centers.

% Indeed, once the estimation $\hat x(t)$ has been computed, the
% condition $\|z(t) - H \hat x(t)\|_{\infty}>\Gamma$ can be checked by
% each monitor without any further communication. A related example is
% presented in Section \ref{sec:numerical}.

\begin{remark}[Statistical detection]
  A different strategy for the detection of false data relies on
  statistical techniques, e.g., see \cite{AA-AGS:04}. In the interest
  of brevity, we do not consider these methods, and we only remark
  that, once the estimation residual has been computed by each
  monitor, the implementation of a (distributed) statistical
  procedure, such as, for instance, the (distributed) $\chi^2$-Test,
  is a straightforward task.
\end{remark}

\section{A finite-memory estimation technique}\label{sec:approximate}
The procedure described in Algorithm \ref{algo:pseudoinverse} allows
each agent to compute an optimal estimate of the whole network state
in finite time. In this section, we allow each agent to handle only
local, i.e., of small dimension, vectors, and we develop a procedure
to recover an estimate of only a certain subnetwork. We envision that
the knowledge of only a subnetwork may be sufficient to implement
distributed estimation and control strategies.

We start by introducing the necessary notation. Let the measurements
matrix $H$ be partitioned into $m^2$, being $m$ the number of monitors
in the network, blocks as
\begin{align}\label{eq:def_H}
  H=\left[ 
    \begin{array}{ccc}
      H_{11} & \cdots & H_{1m}\\
      \vdots & & \vdots \\
      H_{m1} & \cdots & H_{mm}
    \end{array}
  \right], \end{align} where $H_{ij}\in \real^{m_i\times n_i}$ for all
$i,j \in \until{m}$. The above partitioning reflects a division of the
whole network into competence regions: we let each monitor be
responsible for the correct functionality of the subnetwork defined by
its blocks. Additionally, we assume that the union of the different
regions covers the whole network, and that different competence
regions may overlap. Observe that, in most of the practical
situations, the matrix $H$ has a sparse structure, so that many blocks
$H_{ij}$ have only zero entries. We associate an undirected graph
$G_h$ with the matrix $H$, in a way that $G_h$ reflects the
interconnection structure of the blocks $H_{ij}$. To be more precise,
we let $G_h = (V_h , \mathcal{E}_h)$, where $V_h = \until{m}$ denotes
the set of monitors, and where, denoting by $(i,j)$ the undirected
edge from $j$ to $i$, it holds $(i,j) \in \mathcal{E}_h$ if and only
if $\| H_{ij} \| \neq 0$ or $\| H_{ji} \| \neq 0$. Noticed that the
structure of the graph $G_h$, which reflects the sparsity structure of
the measurement matrix, describes also the monitors interconnections.
By using the same partitioning as in \eqref{eq:def_H}, the
Moore-Penrose pseudoinverse of $H$ can be written as
\begin{align}\label{eq:def_Htilde}
  H^\dag= \tilde H =
  \left[ 
    \begin{array}{ccc}
      \tilde{H}_{11} & \cdots & \tilde{H}_{1m}\\
      \vdots & & \vdots \\
      \tilde{H}_{m1} & \cdots & \tilde{H}_{mm}
    \end{array}
  \right],
\end{align}
where $\tilde H_{ij} \in \real^{n_i \times m_i}$. Assume that $H$ has
full row rank,\footnote{The case of a full-column rank matrix is
  treated analogously.} and observe that $H^\dag = H^\transpose (H
H^\transpose)^{-1}$. Consider the equation $z = H x$, and let $H^\dag
z = \hat x = [\hat x_1^\transpose \, \dots \, \hat x_m^\transpose]$,
where, for all $i \in \until{m}$, $\hat x_i \in \real^{n_i}$. We
employ Algorithm \ref{algo:finite_solver} for the computation of the
vector $\hat x$, and we let
\begin{align*}
  \hat x^{(i,h)} = 
  \left[
    \begin{array}{c}
      \hat x_1^{(i,h)}\\
      \vdots\\
      \hat x_m^{(i,h)}
    \end{array}
  \right]
\end{align*}
be the estimate vector of the $i$-th monitor after $h$ iterations of
Algorithm \ref{algo:finite_solver}, i.e., after $h$ executions of the
\emph{while} loop in Algorithm \ref{algo:finite_solver}. In what
follows, we will show that, for a sufficiently sparse matrix $H$, the
error $\| \hat x_i - \hat x_i^{(i,h)} \|$ has an exponential decay
when $h$ increases, so that it becomes negligible before the
termination of Algorithm \ref{algo:finite_solver}, i.e., when $h <
\operatorname{diam}(G_h)$. The main result of this section is next
stated.

\begin{theorem}[Local estimation]\label{thm:local_estimation}
  Let the full-row rank matrix $H$ be partitioned as in
  \eqref{eq:def_H}. Let $[a,b]$, with $a < b$, be the smallest
  interval containing the spectrum of $H H^\transpose$. Then, for $i
  \in \until{m}$ and $h \in \natural$, there exists $C \in \real_{>0}$
  and $q \in (0,1)$ such that
  \begin{align*}
    \| \hat x_i - \hat x_i^{(i,h)} \| \le C q^{\frac{h}{2}+1}.
  \end{align*}
\end{theorem}

Before proving the above result, for the readers convenience, we
recall the following definitions and results. Given an invertible
matrix $M$ of dimension $n$, let us define the \emph{support sets}
\begin{align*}
  \Supp_h(M) = \bigcup_{k=0}^h \{ (i,j) : M^k(i,j) \neq 0 \},
\end{align*}
being $M^k(i,j)$ the $(i,j)$-th entry of $M^k$, and the \emph{decay
  sets}
\begin{align*}
  \Dec_h(M) = (\until{n} \times \until{n}) \setminus \Supp_h(M).
\end{align*}

\begin{theorem}[Decay rate \cite{SD-WFM-PWS:84}]\label{decay_rate}
  Let $M$ be of full row rank, and let $[a,b]$, with $a < b$, be
  the smallest interval containing the spectrum of $M$. % Let $r = b/a$,
  % $q = ( \sqrt{r} - 1 ) / ( \sqrt{r} + 1 )$, and $C = ( 1 + \sqrt{r}
  % )^2 / (2ar)$.
  There exist $C \in \real_{>0}$ and $q \in (0,1)$ such that
\begin{align*}
  \sup \{ | M^\dag (i,j) | : (i,j) \in \Dec_h(M M^\transpose) \} \le C q^{h+1}.
\end{align*}
\end{theorem}

For a graph $G_h$ and two nodes $i$ and $j$, let $\textup{dist}(i,j)$
denote the smallest number of edges in a path from $j$ to $i$ in
$G_h$. The following result will be used to prove Theorem
\ref{thm:local_estimation}. Recall that, for a matrix $M$, we have
$\|M\|_\textup{max} = \max \{ | M(i,j) | \}$.

\begin{lemma}[Decay sets and local
  neighborhood]\label{lemma:dec_neigh}
  Let the matrix $H$ be partitioned as in \eqref{eq:def_H}, and let
  $G_h$ be the graph associated with $H$. For $i,j \in \until{m}$, if
  $\textup{dist}(i,j) = h$, then
  \begin{align*}
    \| H^\dag_{ij}\|_\textup{max} \leq C q^{\frac{h}{2}+1}.
  \end{align*}
\end{lemma}
\begin{pf}
  The proof can be done by simple inspection, and it is omitted here.
\end{pf}

Lemma \ref{lemma:dec_neigh} establishes a relationship between the
decay sets of an invertible matrix and the distance among the vertices
of a graph associated with the same matrix. By using this result, we
are now ready to prove Theorem \ref{thm:local_estimation}.

\begin{pf}[Proof of Theorem \ref{thm:local_estimation}]
  Notice that, after $h$ iterations of Algorithm
  \ref{algo:finite_solver}, the $i$-th monitor has received data from
  the monitors within distance $h$ from $i$, i.e., from the monitors
  $T$ such that, for each $j \in T$, there exists a path of length up
  to $h$ from $j$ to $i$ in the graph associated with $H$. Reorder the
  rows of $H$ such that the $i$-th block come first and the $T$-th
  blocks second. Let $H=[H_1^\transpose \; H_2^\transpose \;
  H_3^\transpose]^\transpose$ be the resulting matrix. Accordingly,
  let $z=[z_1^\transpose \; z_2^\transpose \;
  z_3^\transpose]^\transpose$, and let $x = [x_1^\transpose \;
  x_2^\transpose \; x_3^\transpose]^\transpose$, where $z = H x$.

  Because $H$ has full row rank, we have 
  \begin{align*}
        \left[
          \begin{array}{c}
            H_1\\
            H_2\\
            H_3
          \end{array}
        \right]
        \left[
          \begin{array}{ccc}
            P_{11} & P_{12} & P_{13}\\
            P_{21} & P_{22} & P_{23}\\
            P_{31} & P_{32} & P_{33}
          \end{array}
        \right] = \left[
          \begin{array}{ccc}
            I_1 & 0 & 0\\
            0 & I_2 & 0\\
            0 & 0 & I_3
          \end{array}
        \right],\qquad
        H^\dag = \left[
          \begin{array}{ccc}
            P_{11} & P_{12} & P_{13}\\
            P_{21} & P_{22} & P_{23}\\
            P_{31} & P_{32} & P_{33}
          \end{array}
        \right],
      \end{align*}
      where $I_1$, $I_2$, and $I_3$ are identity matrices of appropriate
      dimension.% , and 
      % \begin{align*}
      %   H^\dag = \left[
      %     \begin{array}{ccc}
      %       P_{11} & P_{12} & P_{13}\\
      %       P_{21} & P_{22} & P_{23}\\
      %       P_{31} & P_{32} & P_{33}
      %     \end{array}
      %   \right].
      % \end{align*}
      For a matrix $M$, let $\textup{col}(M)$ denote the number of
      columns of $M$. Let $T_1 = \until{\textup{col}(P_{11})}$, $T_2 =
      \{1+\textup{col}(P_{11}),\dots,\textup{col}([P_{11} \, P_{12}])
      \}$, and
      \begin{align*}
        T_3 = \{ 1 + \textup{col}([P_{11} \, P_{12}] ) , \dots,
        \textup{col}([P_{11} \, P_{12} \, P_{13}] ) \}.
      \end{align*}
      Let $T_1$, $T_2$, and $T_3$, be, respectively, the indices of
      the columns of $P_{11}$, $P_{12}$, and $P_{13}$. Notice that, by
      construction, if $i \in T_1$ and $j \in T_3$, then
      $\textup{dist} (i,j) > h$. Then, by virtue of Lemma
      \ref{lemma:dec_neigh} and Theorem \ref{decay_rate}, the
      magnitude of each entry of $P_{13}$ is bounded by $\bar C \bar
      q^{\lfloor \frac{h}{2} \rfloor +1}$, for $\bar C, \bar q \in
      \real$.

      Because $H$ has full row rank, from Theorem
      \ref{thm:Algo:Pseudo} we have that
      \begin{align}\label{eq:est}
        \hat x = H^\dag z = \hat{\bar x} + K_1 (H_3 K_1)^\dag (z_3 - H_3 \hat
        x^1),
      \end{align}
      where 
      \begin{align*}
        \hat{\bar x} = [H_1^\transpose \;
        H_2^\transpose]^{\transpose\dag} [z_1^\transpose \;
        z_2^\transpose]^\transpose\;\text{and}\; 
        K_1 = \Basis(
        \Ker([H_1^\transpose \; H_2^\transpose]^\transpose) ).
      \end{align*}
      With the same partitioning as before, let $\hat x = [\hat
      x_1^\transpose \, \hat x_2^\transpose \, \hat
      x_3^\transpose]^\transpose$. In order to prove the theorem, we
      need to show that there exists $C\in \real_{>0}$ and $q \in
      (0,1)$ such that
    \begin{align*}
      \| \hat x_1 - \hat{\bar{x}}_1 \| \le C
      q^{\lfloor \frac{h}{2} \rfloor +1}.
    \end{align*}
    Notice that, for \eqref{eq:est} to hold, the matrix $K_1$ can be
    any basis of $\Ker([H_1^\transpose \;
    H_2^\transpose]^\transpose)$. Hence, let $K_1 = [P_{13}^\transpose
    \; P_{23}^\transpose \; P_{33}^\transpose]$. Because every entry
    of $P_{13}$ decays exponentially, the theorem follows.
\end{pf}

% \begin{figure}
%   \centering
%   \includegraphics[width=.5\columnwidth]{./grid}
%   \caption{A two dimensional power grid with $400$ buses. The network
%     is operated by $16$ control centers, each one responsible for a
%     different subnetwork. Control centers cooperate through the red
%     communication graph.}
%   \label{fig:grid}
% \end{figure}

In Section \ref{sec:simulation_approx} we provide an example to
clarify the exponential decay described in Theorem
\ref{thm:local_estimation}.

\section{An illustrative example}\label{sec:numerical}
The effectiveness of the methods developed in the previous sections is
now shown through some examples.

\subsection{Estimation and detection for the IEEE 118 bus
  system}\label{sec:simulation_118bus}
\begin{figure}[tb]
  \centering \subfigure[]{    
    \includegraphics[width=.46\columnwidth]{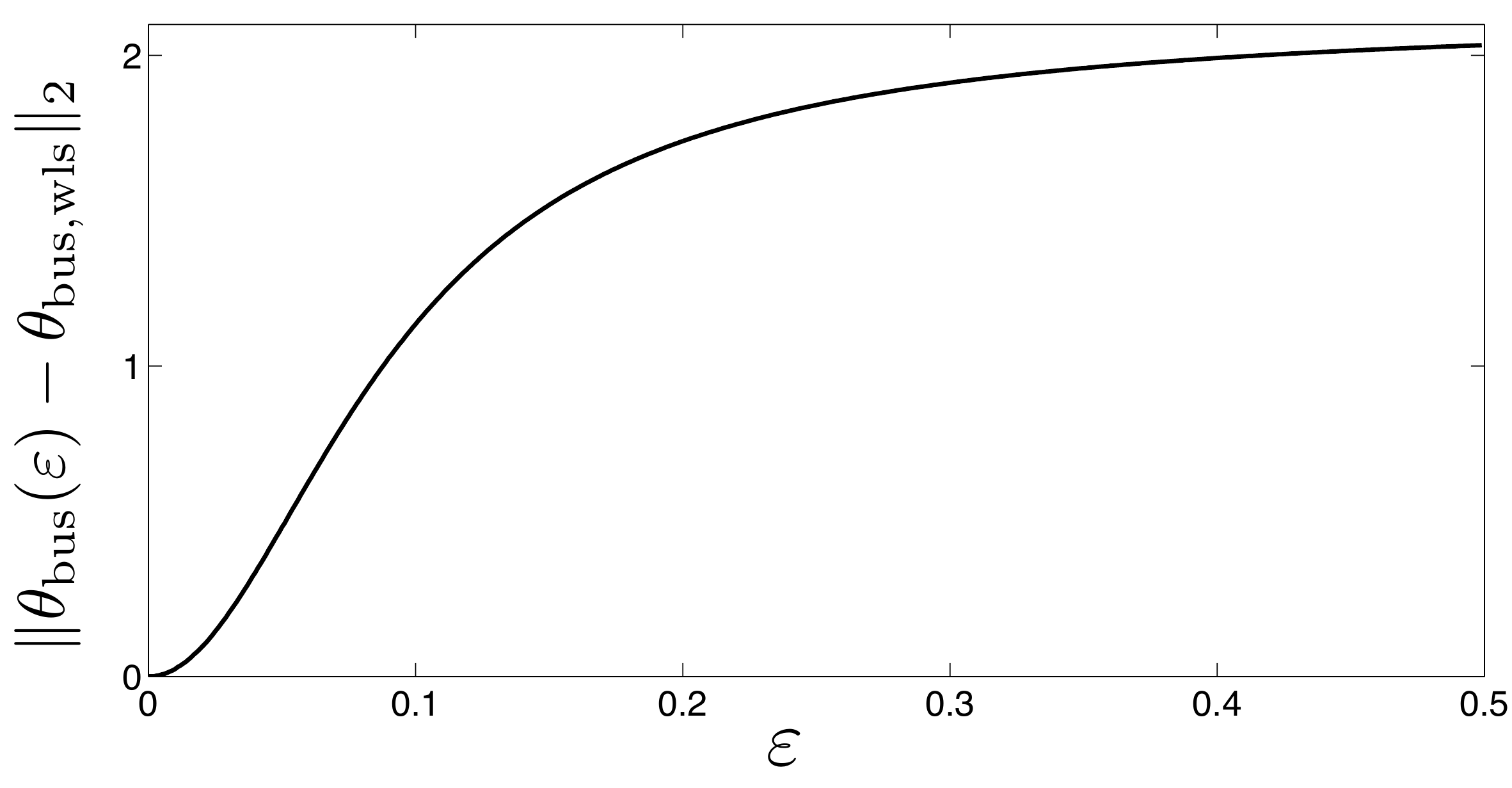}
    \label{fig:error_epsilon}
  } \subfigure[]{
    \includegraphics[width=.5\columnwidth]{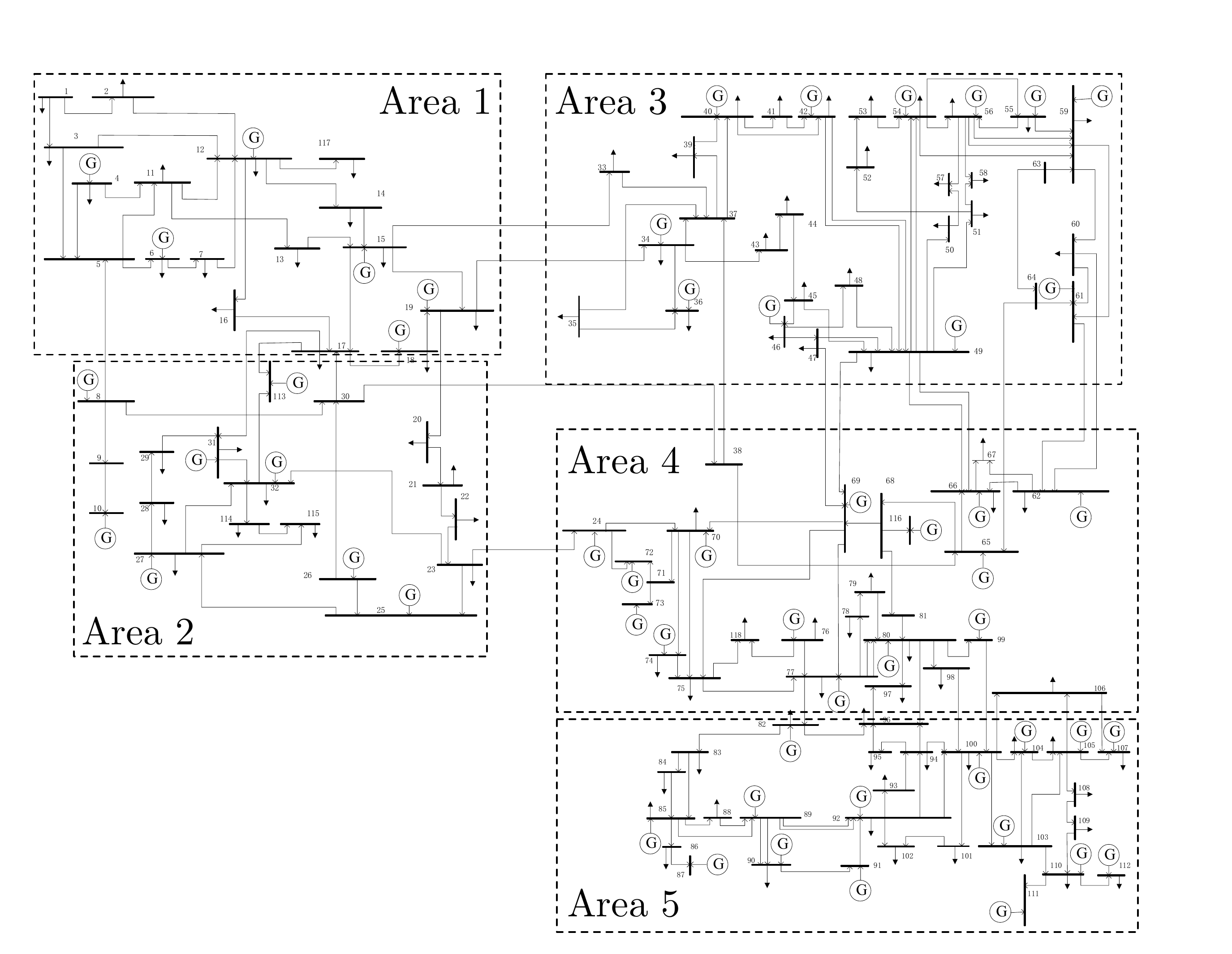}
    \label{fig:network_118_partitioned}
  }
  \caption[Optional caption for list of figures]{In
    Fig. \ref{fig:error_epsilon}, the normalized euclidean norm of the
    error vector $\theta_\textup{bus}(\varepsilon) -
    \theta_\textup{bus,wls}$ is plotted as a function of the parameter
    $\varepsilon$, where $\theta_\textup{bus}(\varepsilon)$ is the
    estimation vector computed according to Theorem
    \ref{existence_limit}, and $\theta_\textup{bus,wls}$ is the
    minimum variance estimate of $\theta_\textup{bus}$. As
    $\varepsilon$ decreases, the vector
    $\theta_\textup{bus}(\varepsilon)$ converges to the minimum
    variance estimate $\theta_\textup{bus,wls}$. In
    Fig. \ref{fig:network_118_partitioned}, the IEEE 118 bus system
    has been divided into $5$ areas. Each area is monitored and
    operated by a control center. The control centers cooperate to
    estimate the state and to assess the functionality of the whole
    network.}
\end{figure}

The IEEE 118 bus system represents a portion of the American Electric
Power System as of December, 1962. This test case system, whose
diagram is reported in Fig. \ref{fig:power_network118}, is composed of
118 buses, 186 branches, 54 generators, and 99 loads. The voltage
angles and the power injections at the network buses are assumed to be
related through the linear relation
\begin{align*}
  P_\textup{bus} = H_\textup{bus} \theta_\textup{bus},
\end{align*}
where the matrix $H_\textup{bus}$ depends upon the network
interconnection structure and the network admittance matrix. For the
network in Fig. \ref{fig:power_network118}, let $z = P_\textup{bus} - v$
be the measurements vector, where $\E[v] = 0$ and $\E[ v v^\transpose]
= \sigma^2 I$, $\sigma \in \real$. Then, following the notation in
Theorem \ref{existence_limit}, the minimum variance estimate of
$\theta_\textup{bus}$ can be recovered as
\begin{align*}
  \lim_{\varepsilon \rightarrow 0^+} [H_\textup{bus} \; \varepsilon
  \sigma I]^\dag z.
\end{align*}
In Fig. \ref{fig:error_epsilon} we show that, as $\varepsilon$
decreases, the estimation vector computed according to Theorem
\ref{existence_limit} converges to the minimum variance estimate of
$\theta_\textup{bus}$.
% \begin{figure}
%     \centering
%     \includegraphics[width=.4\columnwidth]{./error_epsilon}
%     \caption{The normalized euclidean norm of the error vector
%       $\theta_\textup{bus}(\varepsilon) - \theta_\textup{bus,wls}$ is
%       plotted as a function of the parameter $\varepsilon$, where
%       $\theta_\textup{bus}(\varepsilon)$ is the estimation vector
%       computed according to Theorem \ref{existence_limit}, and
%       $\theta_\textup{bus,wls}$ is the minimum variance estimate of
%       $\theta_\textup{bus}$. As $\varepsilon$ decreases, the vector
%       $\theta_\textup{bus}(\varepsilon)$ converges to the minimum
%       variance estimate $\theta_\textup{bus,wls}$.}
%     \label{fig:error_epsilon}
% \end{figure}

% \begin{figure}
%     \centering
%     \includegraphics[width=.6\columnwidth]{./network_118_partitioned}
%     \caption{The IEEE 118 bus system has been divided into $5$
%       areas. Each area is monitored and operated by a control
%       center. The control centers cooperate to estimate the state and
%       to assess the functionality of the whole network.}
%     \label{fig:network_118_partitioned}
% \end{figure}

\begin{figure}[tb]
  \centering \subfigure[]{
    \includegraphics[width=.45\columnwidth]{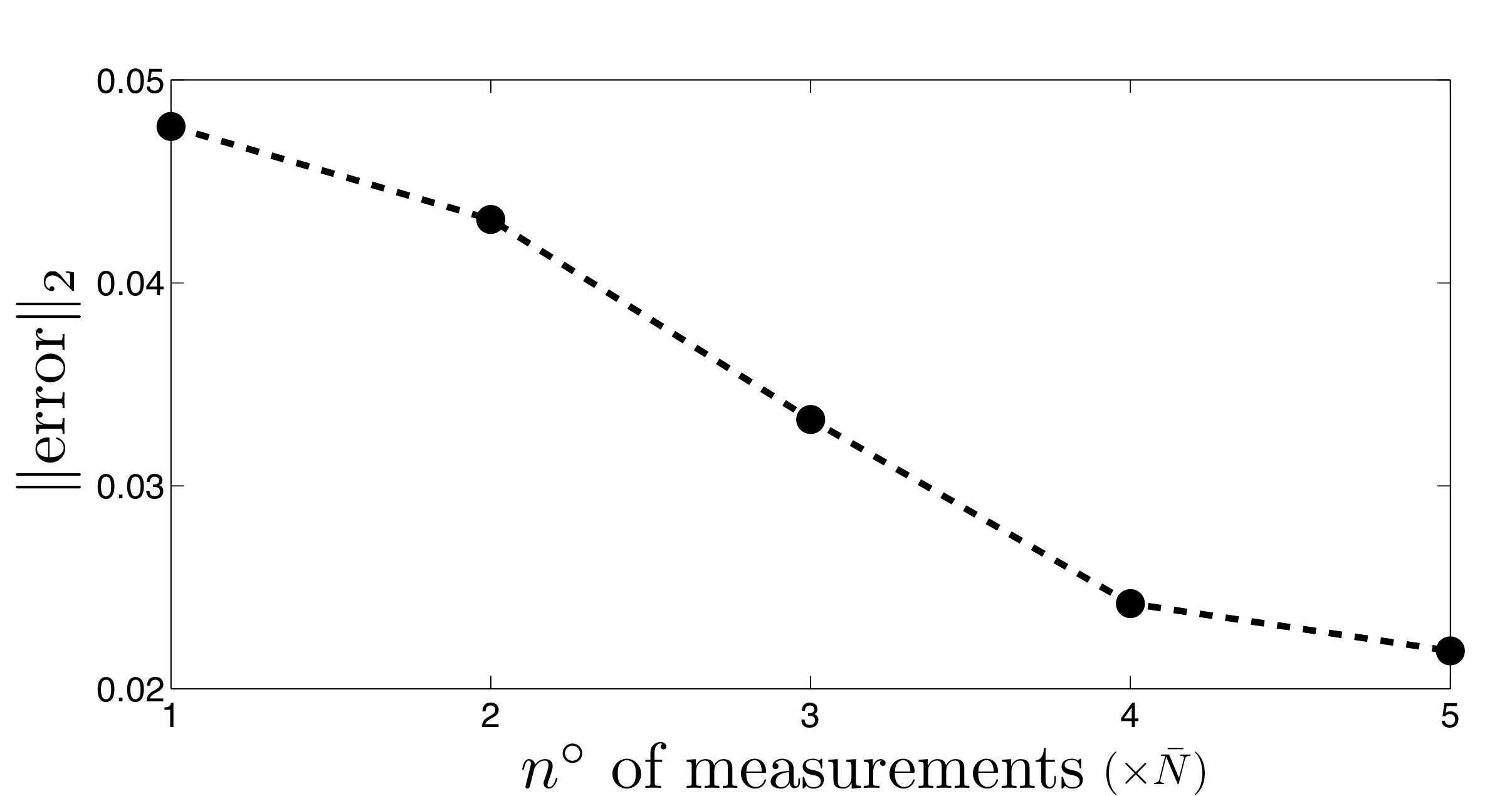}
    \label{fig:increase_measurements}
  } \subfigure[]{
    \includegraphics[width=.49\columnwidth]{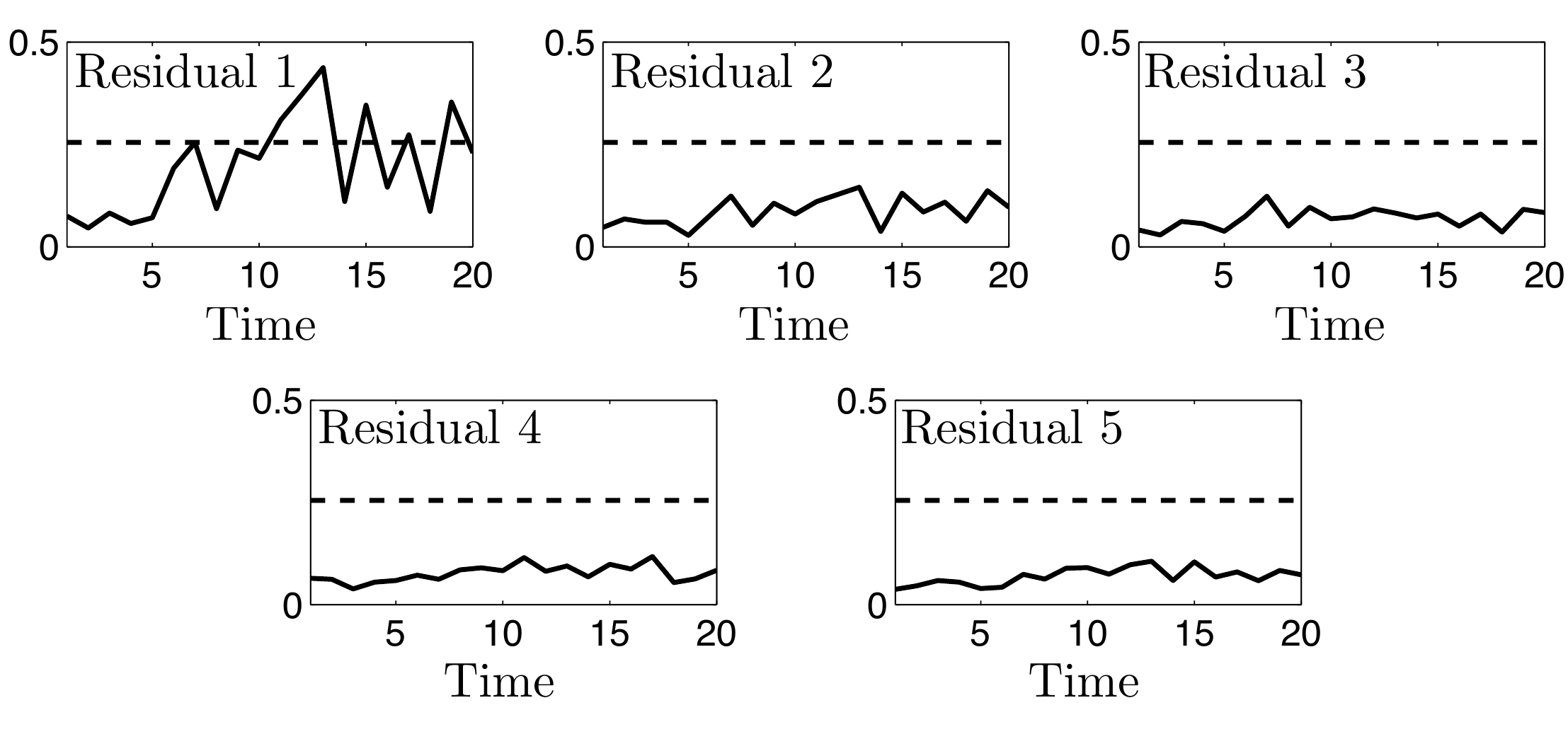}
    \label{fig:detection}
  }
  \caption[Optional caption for list of figures]{For a fixed value of
    $\varepsilon$, Fig. \ref{fig:increase_measurements} shows the
    average (over 100 tests) of the norm of the error (with respect to
    the network state) of the estimate obtained by means of Algorithm
    \ref{algo:pseudoinverse}. The estimation error decreases with the
    number of measurements. Because of the presence of several control
    centers, our distributed algorithm processes more measurements (up
    to $5 \bar N$) while maintaining the same (or smaller)
    computational complexity of a centralized estimation (with $\bar
    N$ measurements). Fig. \ref{fig:detection} shows the residual
    functions computed by the $5$ control centers. Since the first
    residual is greater than the threshold value, the presence of
    false data is correctly detected by the first control center. A
    form of regional identification is possible by simple identifying
    the residuals above the security threshold.}
\end{figure}

In order to demonstrate the advantage of our decentralized estimation
algorithm, we assume the presence of $5$ control centers in the
network of Fig. \ref{fig:power_network118}, each one responsible for a
subpart of the entire network. The situation is depicted in
Fig. \ref{fig:network_118_partitioned}. Assume that each control
center measures the real power injected at the buses in its area, and
let $z_i = P_\textup{bus,i} - v_i$, with $\E[v_i] = 0$ and $\E[ v_i
v_i^\transpose] = \sigma_i^2 I$, be the measurements vector of the
$i$-th area. Finally, assume that the $i$-th control center knows the
matrix $H_\textup{bus,i}$ such that $z_i = H_\textup{bus,i}
\theta_\textup{bus} + v_i$. Then, as discussed in Section
\ref{sec:solver}, the control centers can compute an optimal estimate
of $\theta_\textup{bus}$ by means of Algorithm
\ref{algo:pseudoinverse} or \ref{algo:finite_solver}. Let $n_i$ be the
number of measurements of the $i$-th area, and let $N = \sum_{i=1}^5
n_i$. Notice that, with respect to a centralized computation of the
minimum variance estimate of the state vector, our estimation
procedure obtains the same estimation accuracy while requiring a
smaller computation burden and memory requirement. Indeed, the $i$-th
monitor uses $n_i$ measurements instead of $N$. Let $\bar N$ be the
maximum number of measurements that, due to hardware or numerical
contraints, a control center can efficiently handle for the state
estimation problem. In Fig. \ref{fig:increase_measurements}, we
increase the number of measurements taken by a control center, so that
$n_i \le \bar N$, and we show how the accuracy of the state estimate
increases with respect to a single control center with $\bar N$
measurements.
% \begin{figure}
%     \centering
%     \includegraphics[width=.4\columnwidth]{./increase_measurements}
%     \caption{For a fixed value of $\varepsilon$, the plot shows the
%       norm of the error with respect to the true state vector of the
%       estimate obtained by means of Algorithm \ref{algo:pseudoinverse}
%       (solid) and of the centralized minimum variance estimation with
%       $\bar N$ measurements (dashed). The $x$ axis denotes the number
%       of measurements used for the distributed estimation. Because of
%       the presence of several control centers, the distributed
%       algorithm is more accurate while maintaining the same (or
%       smaller) computational complexity of the centralized
%       estimation.}
%     \label{fig:increase_measurements}
% \end{figure}

To conclude this section, we consider a security application, in which
the control centers aim at detecting the presence of false data among
the network measurements via distributed computation. For this
example, we assume that each control center mesures the real power
injection as well the current magnitude at some of the buses of its
area. By doing so, a sufficient redundancy in the measurements is
obtained for the detection to be feasible \cite{AA-AGS:04}. Suppose
that the measurements of the power injection at the first bus of the
first area is corrupted by a malignant agent. To be more precise, let
the measurements vector of the first area be $\bar z_i = z_i + e_1
w_i$, where $e_1$ is the first canonical vector, and $w_i$ is a random
variable. For the simulation we choose $w_i$ to be uniformly
distributed in the interval $[0,w_\textup{max}]$, where
$w_\textup{max}$ corresponds approximately to the $10\%$ of the
nominal real injection value. In order to detect the presence of false
data among the measurements, the control centers implement Algorithm
\ref{algo:detection}, where, being $H$ the measurements matrix, and
$\sigma$, $\Sigma$ the noise standard deviation and covariance matrix,
the threshold value $\Gamma$ is chosen as $2 \sigma \| I - H
(H^\transpose \Sigma^{-1} H)^{-1} H^\transpose \Sigma^{-1}
\|_{\infty}$.\footnote{For a Gaussian distribution with mean $\mu$ and
  variance $\sigma^2$, about $95\%$ of the realizations are contained
  in $[\mu - 2 \sigma, \mu + 2 \sigma]$.} The residual functions
$\|z_i - H \hat x\|_{\infty}$ are reported in
Fig. \ref{fig:detection}. Observe that, since the first residual is
greater than the threshold $\Gamma$, the control centers successfully
detect the false data. Regarding the identification of the corrupted
measurements, we remark that a regional identification may be possible
by simply analyzing the residual functions. In this example, for
instance, since the residuals $\{2,\dots,5\}$ are below the threshold
value, the corrupted data is likely to be among the measurements of
the first area. This important aspect is left as the subject of future
research.

% \begin{figure}
%     \centering
%     \includegraphics[width=.45\columnwidth]{./detection}
%     \caption{Residual functions computed by the $5$ control
%       centers. Since the first residual is greater than the threshold
%       value, the presence of false data is correctly detected by
%       the first control center. A form of regional identification is
%       possible by simple identifying the residuals above the security
%       threshold.}
%     \label{fig:detection}
% \end{figure}

\subsection{Scalability property of our finite-memory estimation
  technique}\label{sec:simulation_approx}
\begin{figure}[tb]
  \centering \subfigure[]{
    \includegraphics[width=.42\columnwidth]{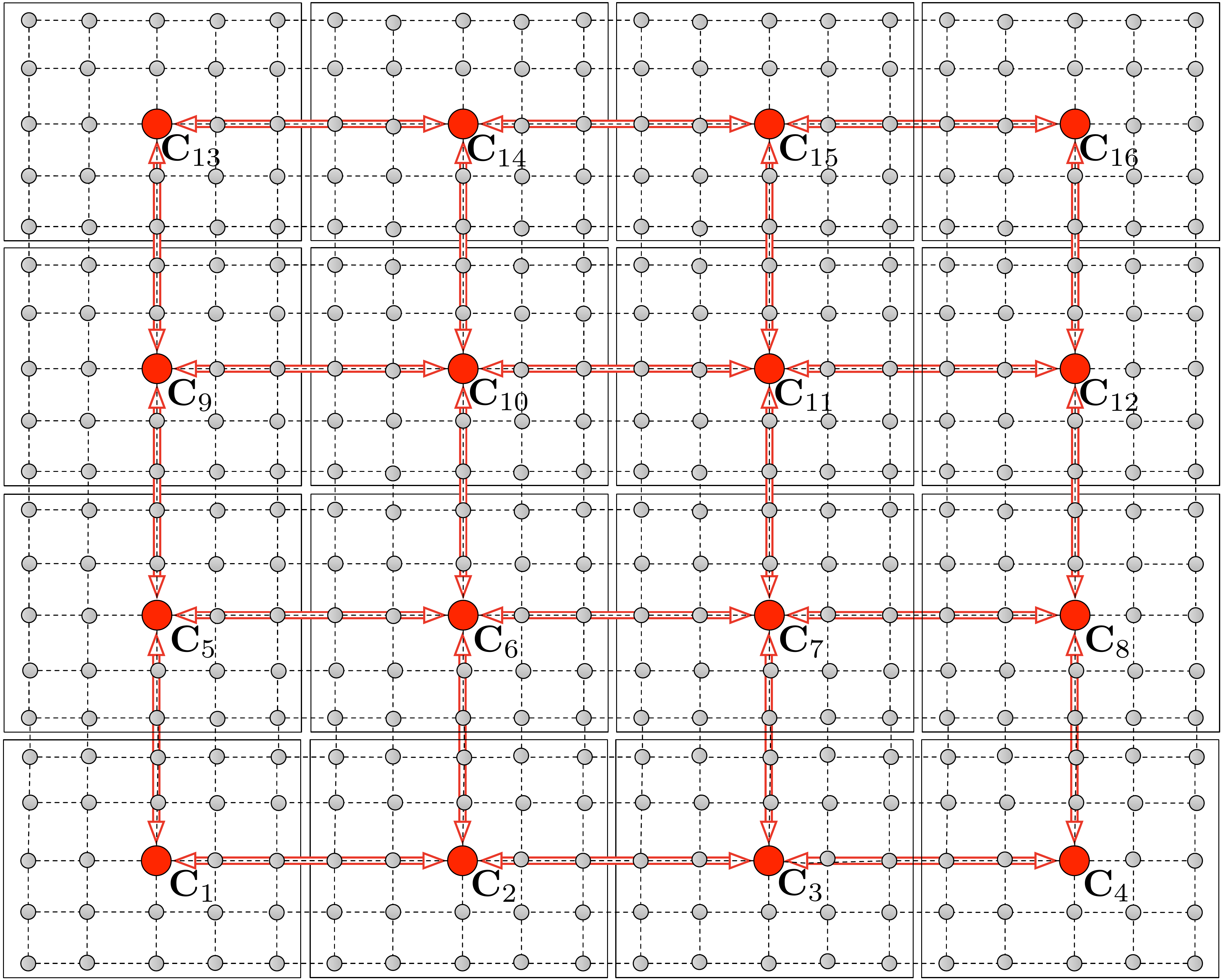}
    \label{fig:grid}
  } \subfigure[]{
    \includegraphics[width=.54\columnwidth]{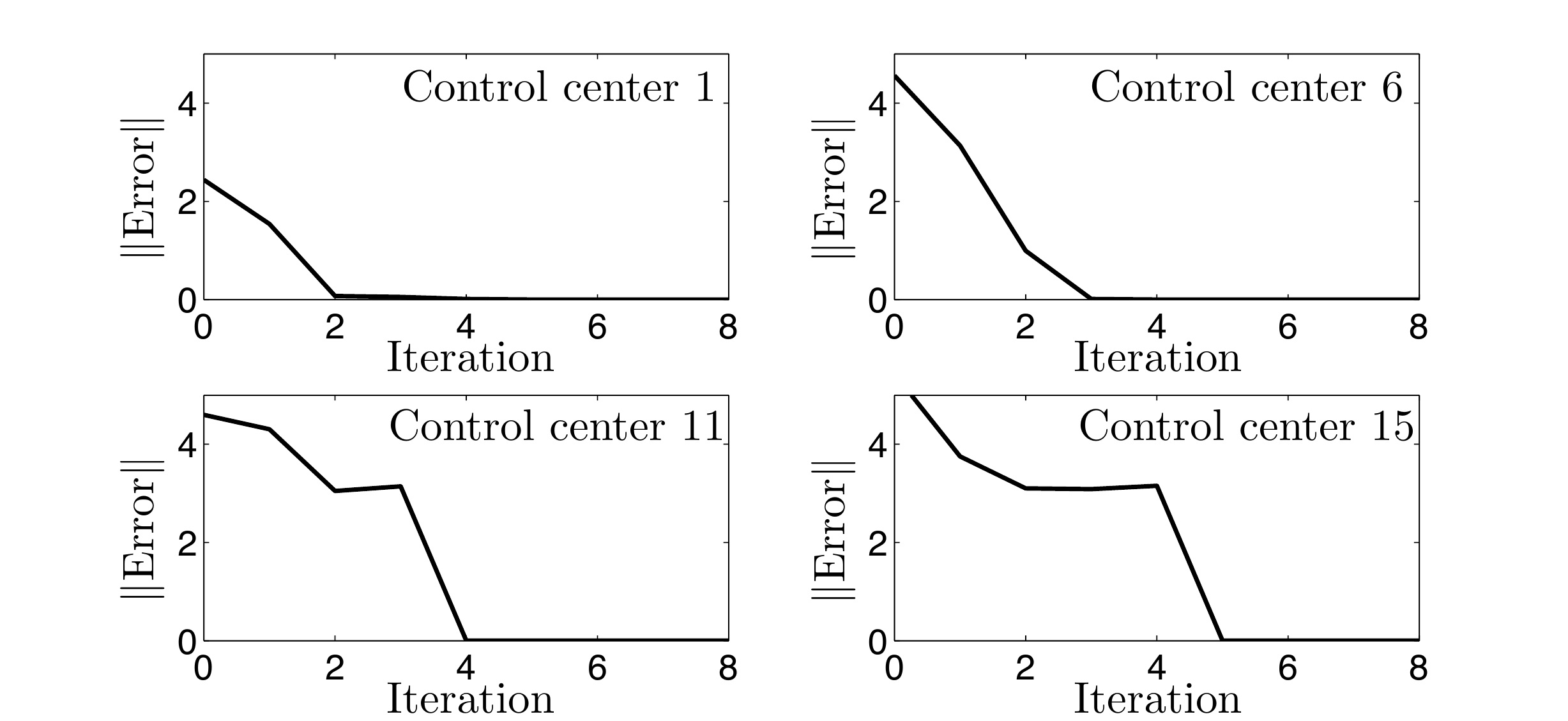}
    \label{fig:simulation_finite}
  }
  \caption[Optional caption for list of figures]{In
    Fig. \ref{fig:grid}, a two dimensional power grid with $400$
    buses. The network is operated by $16$ control centers, each one
    responsible for a different subnetwork. Control centers cooperate
    through the red communication
    graph. Fig. \ref{fig:simulation_finite} shows the norm of the
    estimation error of the local subnetwork as a function of the
    number of iterations of Algorithm \ref{algo:finite_solver}. The
    considered monitors are $C_1$ ,$C_6$, $C_{11}$, and $C_{15}$. As
    predicted by Theorem \ref{thm:local_estimation}, the local
    estimation error becomes negligible before the termination of the
    algorithm.}
\end{figure}

Consider an electrical network with $(ab)^2$ buses, where $a, b \in
\natural$. Let the buses interconnection structure be a two
dimensional lattice, and let $G$ be the graph whose vertices are the
$(ab)^2$ buses, and whose edges are the network branches. Let $G$ be
partitioned into $b^2$ identical blocks containing $a^2$ vertices
each, and assume the presence of $b^2$ control centers, each one
responsible for a different network part. We assume the control
centers to be interconnected through an undirected graph. In
particular, being $V_i$ the set of buses assigned to the control
center $\textup{C}_i$, we let the control centers $\textup{C}_i$ and
$\textup{C}_j$ be connected if there exists a network branch linking a
bus in $V_i$ to a bus in $V_j$. An example with $b = 4$ and $a = 5$ is
in Fig. \ref{fig:grid}. In order to show the effectiveness of our
approximation procedure, suppose that each control center
$\textup{C}_i$ aims at estimating the vector of the voltage angles at
the buses in its region. We assume also that the control centers
cooperate, and that each of them receives the measurements of the real
power injected at only the buses in its region. Algorithm
\ref{algo:finite_solver} is implemented by the control centers to
solve the estimation problem. In Fig. \ref{fig:simulation_finite} we
report the estimation error during the iterations of the
algorithm. Notice that, as predicted by Theorem
\ref{thm:local_estimation}, each leader possess a good estimate of the
state of its region before the termination of the algorithm.

\section{Conclusion}\label{sec:conclusion}
Two distributed algorithms for network control centers to compute the
minimum variance estimate of the network state given noisy
measurements have been proposed. The two methods differ in the mode of
cooperation of the control centers: the first method implements an
incremental mode of cooperation, while the second uses a diffusive
interaction. Both methods converge in finite time, which we
characterize, and they require only local measurements and model
knowledge to be implemented. Additionally, an asynchronous and
scalable implementation of our diffusive estimation method has been
described, and its efficiency has been shown through a rigorous
analysis and through a practical example. Based on these estimation
methods, an algorithm to detect cyber-attacks against the network
measurements has also been developed, and its detection performance
has been characterized.

\renewcommand{\theequation}{A-\arabic{equation}}
% redefine the command that creates the equation no.
\setcounter{equation}{0}  % reset counter 
\section*{APPENDIX}  % use *-form to suppress numbering
%\section{Appendix}\label{appendix}
\subsection{Proof of Theorem
  \ref{thm:Algo:Pseudo}}\label{pf_algopseudo}

\begin{pf}
  % For $1 \le k\le m$, let
%   \begin{align*}
%     H^k=\left[
%       \begin{array}{c}
%         H_1\\
%         H_2\\
%         \vdots\\
%         H_k
%       \end{array}
%     \right],\quad
%     z^k=\left[
%       \begin{array}{c}
%         z_1\\
%         z_2\\
%         \vdots\\
%         z_k
%       \end{array}
%     \right].
% \end{align*}
% The next result shows that the vector $\hat x_i$ and the matrix $K_i$
% in Algorithm \ref{algo:pseudoinverse} are such that $y^i=H^i \hat x_i$
% and $\hat x_i \perp \Ker(H^i)$ for all $1 \le i \le m$, so that, at
% each step $k$, Algorithm \ref{algo:pseudoinverse} computes the minimum
% norm solution of the system of linear equations $z^k = H^k x$.
  Let $H^i=[H_1^\transpose \enspace \cdots \enspace
  H_i^\transpose]^\transpose$, $z^i=[z_1^\transpose \enspace \cdots
  \enspace z_i^\transpose]^\transpose$. We show by induction that
  $z^i=H^i \hat x_i$, $K_i = \Basis(\Ker(H^i))$, and $\hat x_i \perp
  \Ker(H^i)$. Note that the statements are trivially verified for
  $i=1$. Suppose that they are verified up to $i$, then we need to
  show that $K_{i+1} = \Basis(\Ker(H^{i+1}))$, $\hat x_{i+1} \perp
  \Ker(H^{i+1})$, and $z^{i+1}=H^{i+1} \hat x_{i+1}$. 

  We start by proving that $K_{i+1} = \Basis(\Ker(H^{i+1}))$. Observe
  that $\Ker(K_i)=0$ for all $i$, and that
 \begin{align}\label{eq:kernel}
    \Ker(H_{i+1}K_i) = \{v  :  K_i v \in \Ker(H_{i+1})\}.
 \end{align}
 Hence,
% It follows that $\Image(K_{i+1}) = \Image(K_i \Ker(H_{i+1}K_i))
% \subseteq \Ker(H_{i+1})$, and, since trivially $K_i w \in
 %\Image(K_i)$ for every vector $w$ of appropriate dimension, we have
  \begin{align*}
    \Image(K_{i+1}) &= \Image(K_i \Ker(H_{i+1}K_i)) = \Image(K_i) \cap \Ker(H_{i+1}) = \Ker(H^{i})
    \cap \Ker(H_{i+1}) = \Ker(H^{i+1}).
  \end{align*}

  We now show that $x_{i+1} \perp \Ker(H^{i+1})$, which is equivalent to
  \begin{align*}
    \hat x_{i+1}= (\hat x_i + K_i(H_{i+1}K_i)^\dag(z_{i+1}-H_{i+1}
    \hat x_i) ) \in \Ker(H^{i+1})^{\perp}.
  \end{align*}
  Note that
  \begin{align*}
    \Ker(H^{i+1}) \subseteq \Ker(H^i) \enspace \Leftrightarrow
    \enspace \Ker(H^{i+1})^\perp \supseteq \Ker(H^i)^\perp.
  \end{align*}
  By the induction hypothesis we have $\hat x_i \in \Ker(H^i)^\perp$,
  and hence $\hat x_i \in \Ker(H^{i+1})^\perp$.  
Therefore, we need to show that
  \begin{align*}
    K_i(H_{i+1}K_i)^\dag(z_{i+1}-H_{i+1} \hat x_i) \in
    \Ker(H^{i+1})^{\perp}.
  \end{align*}
  Let $w = (H_{i+1}K_i)^\dag(z_{i+1}-H_{i+1} \hat x_i)$, and notice
  that $w \in \Ker(H_{i+1} K_i)^\perp$ due to the properties of the
  pseudoinverse operation. Suppose that $K_i w \not\in
  \Ker(H_{i+1})^\perp$. Since $\Ker(K_i) = \{0\}$, the vector $w$ can
  be written as $w = w_1 + w_2$, where $K_i w_1 \in
  \Ker(H_{i+1})^\perp$ and $K_i w_2 = K_i w - K_i w_1 \neq 0$, $K_i
  w_2 \in \Ker(H_{i+1})$. Then, it holds $H_{i+1} K_i w_2 = 0$, and
  hence $w_2 \in \Ker(H_{i+1} K_i )$, which contradicts the hypothesis
  $w \in \Ker(H_{i+1} K_i)^\perp$. Finally $K_i w \in
  \Ker(H_{i+1})^\perp \subseteq \Ker(H^{i+1})^\perp$. 

  We now show that $z^{i+1}=H^{i+1} \hat x_{i+1}$. Because of the
  consistency of the system of linear equations, and because $z^{i} =
  H^{i} \hat x_i$ by the induction hypothesis, there exists a vector
  $v_i \in \Ker(H^{i})=\Image(K_i)$ such that $z^{i+1} = H^{i+1} (\hat
  x_i + v_i)$, and hence that $z_{i+1} = H_{i+1} (\hat x_i + v_i)$. We
  conclude that $(z_{i+1}-H_{i+1} \hat x_i) \in \Image (H_{i+1} K_i)$,
  and finally that $z^{i+1} = H^{i+1} \hat x_{i+1}$.

  % Notice that $(H_{i+1}K_i)^\dag(z_{i+1}-H_{i+1} \hat x_i) = w \perp
  % \Ker(H_{i+1} K_i)$. Hence, from equation \eqref{eq:kernel}, we have
  % $K_i w \perp \Ker(H_{i+1})$, and, finally, $K_i w \perp
  % \Ker(H^{i+1})$ and $\hat x_{i+1} \perp \Ker(H^{i+1})$. Because of
  % the consistency of the system of linear equations, and because
  % $z^{i} = H^{i} \hat x_i$ by the induction hypothesis, there exists a
  % vector $v_i \in \Ker(H^{i})=\Image(K_i)$ such that $z^{i+1} =
  % H^{i+1} (\hat x_i + v_i)$, and hence that $z_{i+1} = H_{i+1} (\hat
  % x_i + v_i)$. We conclude that $(z_{i+1}-H_{i+1} \hat x_i) \in \Image
  % (H_{i+1} K_i)$, and finally that $z^{i+1} = H^{i+1} \hat x_{i+1}$.

% Therefore, we need to show that
%   \begin{align*}
%     K_i(H_{i+1}K_i)^\dag(z_{i+1}-H_{i+1} \hat x_i) \in
%     \Ker(H^{i+1})^{\perp},
%   \end{align*}
%   which is implied by
%   \begin{align*}
%     H_{i+1} K_i(H_{i+1}K_i)^\dag (z_{i+1}-H_{i+1} \hat x_i) \neq 0,
%   \end{align*}
%   or by
%   \begin{align*}
%     (z_{i+1}-H_{i+1} \hat x_i) \in \Image (H_{i+1} K_i).
%   \end{align*}
%   Because of the consistency of the system of linear equations, and
%   because $z^{i} = H^{i} \hat x_i$ by the induction hypothesis, there
%   exists a vector $v_i \in \Ker(H^{i})=\Image(K_i)$ such that $z^{i+1}
%   = H^{i+1} (\hat x_i + v_i)$, and hence that $z_{i+1} = H_{i+1} (\hat
%   x_i + v_i)$. We conclude that $(z_{i+1}-H_{i+1} \hat x_i) \in \Image
%   (H_{i+1} K_i)$, and finally that $z^{i+1} = H^{i+1} \hat x_{i+1}$
%   and $\hat x_{i+1} \perp \Ker(H^{i+1})$.
\end{pf}

\subsection{Proof of Theorem \ref{existence_limit}}\label{proof_existence_limit}
Before proceeding with the proof of the above theorem, we recall the
following fact in linear algebra.

\begin{lemma}\label{lemma_ker}
  Let $H \in \mathbb{R}^{n\times m}$. Then
  $\Ker((H^\dag)^\transpose)=\Ker(H)$.
\end{lemma}
\begin{pf}
  We first show that $\Ker((H^\dag)^\transpose)\subseteq
  \Ker(H)$. Recall from \cite{DSB:09} that
  $H=HH^\transpose(H^\dag)^\transpose$. Let $x$ be such that
  $(H^\dag)^\transpose x = 0$, then $H x =
  HH^\transpose(H^\dag)^\transpose x =0$, so that
  $\Ker((H^\dag)^\transpose)\subseteq \Ker(H)$. We now show that
  $\Ker(H) \subseteq \Ker((H^\dag)^\transpose)$. Recall that
  $(H^\dag)^\transpose=(H^\transpose)^\dag=(HH^\transpose)^\dag
  H$. Let $x$ be such that $H x =0$, then $(H^\dag)^\transpose x =
  (HH^\transpose)^\dag H x =0$, so that $\Ker(H) \subseteq
  \Ker((H^\dag)^\transpose)$, which concludes the proof.
\end{pf}
We are now ready to prove Theorem \ref{existence_limit}.
\begin{pf}
  The first property follows directly from \cite{DSB:09} (cfr. page
  427). To show the second property, observe that $C^\dag =
  \frac{1}{\varepsilon} ((I - HH^\dag)B)^\dag$, so that
\begin{align*}
  \lim_{\varepsilon \rightarrow 0^+} \varepsilon D=0.
\end{align*}
For the theorem to hold, we need to verify that
\begin{align*}
  H^\dag - H^\dag B ((I - HH^\dag)B)^\dag = (H^\transpose \Sigma^{-1}
  H)^{-1}H^\transpose\Sigma^{-1},
\end{align*}
or, equivalently, that
\begin{align}\label{first}
  \begin{split}
    &\left(H^\dag -H^\dag B ((I-HH^\dag)B)^\dag\right) HH^\dag =
    (H^\transpose \Sigma^{-1} H)^{-1}H^\transpose \Sigma^{-1} HH^\dag,
  \end{split}
\end{align}
and
\begin{align}\label{second}
  \begin{split}
    &\left(H^\dag -H^\dag B ((I-HH^\dag)B)^\dag\right) (I - HH^\dag
    )= (H^\transpose \Sigma^{-1} H)^{-1}H^\transpose \Sigma^{-1}
    (I - HH^\dag ).
  \end{split}
\end{align}

Consider equation \eqref{first}. After simple manipulation, we have
\begin{align*}
  H^\dag - H^\dag B ((I-HH^\dag)B)^\dag HH^\dag &=H^\dag,
 \end{align*}
so that we need to show only that
\begin{align*}
  H^\dag B ((I-HH^\dag)B)^\dag HH^\dag &= 0.
\end{align*}
Recall that for a matrix $W$ it holds $W^\dag = (W^\transpose W)^\dag
W^\transpose$. Then the term $((I-HH^\dag)B)^\dag HH^\dag$ equals
\begin{align*}
  \left( ((I-HH^\dag)B)^\transpose ((I-HH^\dag)B) \right)^\dag
  B^\transpose (I-HH^\dag) HH^\dag = 0,
\end{align*}
because $(I-HH^\dag) HH^\dag = 0$. We conclude that equation
\eqref{first} holds. Consider now equation \eqref{second}. Observe
that $HH^\dag (I-HH^\dag)=0$. Because $B$ has full row rank, and
$\Sigma=BB^\transpose$, simple manipulation yields
\begin{align*}
  &-H^\transpose(BB^\transpose)^{-1}HH^\dag B \left[ (I-HH^\dag)B
  \right]^\dag (I-HH^\dag) B= H^\transpose(BB^\transpose)^{-1}
  (I-HH^\dag) B,
\end{align*}
and hence
\begin{align*}
  &H^\transpose(BB^\transpose)^{-1} \left\{ I + HH^\dag B \left[
      (I-HH^\dag)B \right]^\dag \right\} (I-HH^\dag) B=0.
\end{align*}
Since $HH^\dag=I - (I-HH^\dag)$, we obtain
\begin{align*}
  H^\transpose(BB^\transpose)^{-1}B \left[ (I-HH^\dag)B \right]^\dag
  (I-HH^\dag) B =0.
\end{align*}
A sufficient condition for the above equation to be true is
% \begin{align*}
%   A^\transpose(BB^\transpose)^{-1}B \left[ (I-AA^\dag)B \right]^\dag=0.
% \end{align*}
% By taking the transpose of both sides we obtain
\begin{align*}
  \left(\left[ (I-HH^\dag)B \right]^\dag\right)^\transpose
  B^\transpose(BB^\transpose)^{-1} H=0.
\end{align*} 
% Recall from \cite{DSB:09} (cfr. page 399) that 
From Lemma \ref{lemma_ker} we have.
\begin{align*}
  \Ker\left(\left(\left[ (I-AA^\dag)B
      \right]^\dag\right)^\transpose\right)=\Ker((I-AA^\dag)B).
\end{align*}
Since 
\begin{align*}
  (I-HH^\dag)B B^\transpose(BB^\transpose)^{-1} H=(I-HH^\dag)H=0,
\end{align*}
we have that
\begin{align*}
  H^\transpose(BB^\transpose)^{-1}B \left[ (I-HH^\dag)B \right]^\dag
  (I-HH^\dag) B =0,
\end{align*}
and that equation \eqref{second} holds. This concludes the proof.
\end{pf}


\begin{thebibliography}{10}

\bibitem{AA-AGS:04}
A.~Abur and A.~G. Exposito.
\newblock {\em Power System State Estimation: Theory and Implementation}.
\newblock CRC Press, 2004.

\bibitem{DSB:09}
D.~S. Bernstein.
\newblock {\em Matrix Mathematics}.
\newblock Princeton University Press, 2 edition, 2009.

\bibitem{RC-AC-LS-SZ:08a}
R.~Carli, A.~Chiuso, L.~Schenato, and S.~Zampieri.
\newblock Distributed {K}alman filtering based on consensus strategies.
\newblock {\em IEEE Journal on Selected Areas in Communications},
  26(4):622--633, 2008.

\bibitem{FSC-CGL-AHS:08}
F.~S. Cattivelli, C.~G. Lopes, and A.~H. Sayed.
\newblock Diffusion recursive least-squares for distributed estimation over
  adaptive networks.
\newblock {\em IEEE Transactions on Signal Processing}, 56(5):1865--1877, 2008.

\bibitem{FSC-AHS:10}
F.~S. Cattivelli and A.~H. Sayed.
\newblock Diffusion strategies for distributed {K}alman filtering and
  smoothing.
\newblock {\em IEEE Transactions on Automatic Control}, 55(9):2069--2084, 2010.

\bibitem{YC:81}
Y.~Censor.
\newblock Row-action methods for huge and sparse systems and their
  applications.
\newblock {\em SIAM Review}, 23(4):444--466, 1981.

\bibitem{SD-WFM-PWS:84}
S.~Demko, W.~F. Moss, and P.~W. Smith.
\newblock {Decay rates for inverses of band matrices}.
\newblock {\em Mathematics of Computation}, 43(168):491--499, 1984.

\bibitem{DMF-FFW-LM:95}
D.~M. Falcao, F.~F. Wu, and L.~Murphy.
\newblock Parallel and distributed state estimation.
\newblock {\em IEEE Transactions on Power Systems}, 10(2):724--730, 1995.

\bibitem{CDG-GFR:01}
C.~D. Godsil and G.~F. Royle.
\newblock {\em Algebraic Graph Theory}, volume 207 of {\em Graduate Texts in
  Mathematics}.
\newblock Springer, 2001.

\bibitem{GHG-CFvL:89}
G.~H. Golub and C.~F. van Loan.
\newblock {\em Matrix Computations}.
\newblock Johns Hopkins University Press, 2 edition, 1989.

\bibitem{RG-RB-GTH:70}
R.~Gordon, R.~Bender, and G.~T. Herman.
\newblock Algebraic reconstruction techniques ({ART}) for three-dimensional
  electron microscopy and x-ray photography.
\newblock {\em Journal of theoretical Biology}, 29(3):471--481, 1970.

\bibitem{WJ-VV-GTH:07}
W.~Jiang, V.~Vittal, and G.~T. Heydt.
\newblock A distributed state estimator utilizing synchronized phasor
  measurements.
\newblock {\em IEEE Transactions on Power Systems}, 22(2):563--571, 2007.

\bibitem{SK:37}
S.~Kaczmarz.
\newblock {Angen\"aherte Aufl\"osung von Systemen linearer Gleichungen}.
\newblock {\em Bull. Acad. Polon. Sci. Lett. A}, 35:355--357, 1937.

\bibitem{YL-MKR-PN:09}
Y.~Liu, M.~K. Reiter, and P.~Ning.
\newblock False data injection attacks against state estimation in electric
  power grids.
\newblock In {\em ACM Conference on Computer and Communications Security},
  pages 21--32, Chicago, IL, USA, November 2009.

\bibitem{CGL-AHS:07}
C.~G. Lopes and A.~H. Sayed.
\newblock Incremental adaptive strategies over distributed networks.
\newblock {\em IEEE Transactions on Signal Processing}, 55(8):4064--4077, 2007.

\bibitem{CGL-AHS:08}
C.~G. Lopes and A.~H. Sayed.
\newblock Diffusion least-mean squares over adaptive networks: {F}ormulation
  and performance analysis.
\newblock {\em IEEE Transactions on Signal Processing}, 56(7):3122--3136, 2008.

\bibitem{DGL:69}
D.~G. Luenberger.
\newblock {\em Optimization by Vector Space Methods}.
\newblock Wiley, 1969.

\bibitem{JM:07}
J.~Meserve.
\newblock Sources: Staged cyber attack reveals vulnerability in power grid.
\newblock http://cnn.com, September 26, 2007.

\bibitem{AM:99}
A.~Monticelli.
\newblock {\em State Estimation in Electric Power Systems: A Generalized
  Approach}.
\newblock Springer, 1999.

\bibitem{N:04}
{NERC}.
\newblock {Final Report on the August 14, 2003 Blackout in the United States
  and Canada: Causes and Recommendations}, April 2004.
\newblock {Available at http://www.nerc.com/filez/blackout.html}.

\bibitem{FB-RC-AB-FB:10n}
F.~Pasqualetti, R.~Carli, A.~Bicchi, and F.~Bullo.
\newblock Distributed estimation and detection under local information.
\newblock In {\em IFAC Workshop on Distributed Estimation and Control in
  Networked Systems}, pages 263--268, Annecy, France, September 2010.

\bibitem{MGR-RDN:05}
M.~G. Rabbat and R.~D. Nowak.
\newblock Quantized incremental algorithms for distributed optimization.
\newblock {\em IEEE Journal on Selected Areas in Communications},
  23(4):798--808, 2005.

\bibitem{CR-SP-SU:05}
C.~Rakpenthai, S.~Premrudeepreechacharn, S.~Uatrongjit, and N.~R. Watson.
\newblock {Measurement placement for power system state estimation using
  decomposition technique}.
\newblock {\em Electric Power Systems Research}, 75(1):41--49, 2005.

\bibitem{AHS-CGL:07}
A.~H. Sayed and C.~G. Lopes.
\newblock Adaptive processing over distributed networks.
\newblock {\em IEICE Transactions on Fundamentals of Electronics,
  Communications and Computer Sciences}, E90-A(8):1504--1510, 2007.

\bibitem{IS-AR-GG:07}
I.~Schizas, A.~Ribeiro, and G.~Giannakis.
\newblock Consensus in ad hoc {WSN}s with noisy links - {P}art {I}: Distributed
  estimation of deterministic signals.
\newblock {\em IEEE Transactions on Signal Processing}, 56(1):350--364, 2007.

\bibitem{IDS-GM-GBG:09}
I.~D. Schizas, G.~Mateos, and G.~B. Giannakis.
\newblock Distributed {LMS} for consensus-based in-network adaptive processing.
\newblock {\em IEEE Transactions on Signal Processing}, 57(6):2365--2382, 2009.

\bibitem{FCS-JW:70}
F.~C. Schweppe and J.~Wildes.
\newblock Power system static-state estimation, {Part I: E}xact model.
\newblock {\em IEEE Transactions on Power Apparatus and Systems},
  89(1):120--125, 1970.

\bibitem{FCS-JW:70-bis}
F.~C. Schweppe and J.~Wildes.
\newblock Power system static-state estimation, {Part II: A}pproximate model.
\newblock {\em IEEE Transactions on Power Apparatus and Systems},
  89(1):125--130, 1970.

\bibitem{FCS-JW:70-bisbis}
F.~C. Schweppe and J.~Wildes.
\newblock Power system static-state estimation, {Part III: I}mplementation.
\newblock {\em IEEE Transactions on Power Apparatus and Systems},
  89(1):130--135, 1970.

\bibitem{MS-YW:03}
M.~Shahidehpour and Y.~Wang.
\newblock {\em Communication and Control in Electric Power Systems:
  Applicationsc of Parallel and Distributed Processing}.
\newblock Wiley-IEEE Press, 2003.

\bibitem{SSS-MSS-DMS:09}
S.~S. Stankovic, M.~S. Stankovic, and D.~M. Stipanovic.
\newblock Consensus based overlapping decentralized estimation with missing
  observations and communication faults.
\newblock {\em Automatica}, 45(6):1397--1406, 2009.

\bibitem{KT:71}
K.~Tanabe.
\newblock Projection method for solving a singular system of linear equations
  and its applications.
\newblock {\em Numerische Mathematik}, 17(3):203--214, 1971.

\bibitem{LX-YM-BS:10}
L.~Xie, Y.~Mo, and B.~Sinopoli.
\newblock False data injection attacks in electricity markets.
\newblock In {\em IEEE Int. Conf. on Smart Grid Communications}, pages
  226--231, Gaithersburg, MD, USA, October 2010.

\bibitem{LZ-AA:05}
L.~Zhao and A.~Abur.
\newblock Multi area state estimation using synchronized phasor measurements.
\newblock {\em IEEE Transactions on Power Systems}, 20(2):611--617, 2005.

\end{thebibliography}
\end{document}